\documentclass{amsart}
\usepackage{amsmath,amsfonts,latexsym,amssymb,amsthm}
\usepackage[mathscr]{eucal}

\newtheorem{theorem}{Theorem}[section]
\newtheorem{corollary}[theorem]{Corollary}

\newtheorem{lemma}[theorem]{Lemma}
\newtheorem*{main}{Theorem~\ref{newthm}}
\newtheorem*{question}{Question}

\theoremstyle{definition}
\newtheorem{definition}[theorem]{Definition}

\newtheorem{remark}[theorem]{Remark}

\theoremstyle{definition}

\newtheorem{clm}{Claim}

\theoremstyle{definition}
\newtheorem*{clam}{Claim}

\theoremstyle{remark}
\newtheorem*{rmk}{Comment}

\newcommand{\be}{\begin{enumerate}}
\newcommand{\ee}{\end{enumerate}}
\newcommand{\bi}{\begin{itemize}}
\newcommand{\ei}{\end{itemize}}
\newcommand{\seq}[1]{\langle#1\rangle}

\newcommand{\R}{\ensuremath{\mathbb{R}} }
\newcommand{\n}{\ensuremath{{n}}}
\newcommand{\m}{\ensuremath{{r}}}
\newcommand{\ul}{\underline}
\renewcommand{\H}{\mathcal{H}}
\newcommand{\en}{\text{$\in$}}
\newcommand{\lh}{\textup{lh}}

\def\F{\mathfrak{F}}
\def\Eq{{\underline{E}}}
\def\f{\mathfrak{f}}
\def\g{\mathfrak{g}}
\def\h{\mathfrak{h}}

\def\pow{{\mathcal{P}}}
\def\Lng{\mathcal{L}}
\def\G{\mathcal{G}}
\def\l{[\![}
\def\r{]\!]}
\def\Lr{{\ensuremath{L({\R})} }}
\def\Kr{{\ensuremath{K({\R})} }}
\newcommand{\A}{\ensuremath{\mathcal{A}}}
\newcommand{\B}{\ensuremath{\mathcal{B}}}
\newcommand{\fB}{\ensuremath{\mathfrak{B}}}
\newcommand{\fA}{\ensuremath{\mathfrak{A}}}
\newcommand{\T}{\ensuremath{\mathcal{T}}}
\newcommand{\mouse}{\ensuremath{\mathcal{M}} }
\newcommand{\nouse}{\ensuremath{\mathcal{N}} }

\newcommand{\fnouse}{\ensuremath{\mathfrak{N}}}

\newcommand{\overmouse}{\ensuremath{\overline{\mouse}} }
\newcommand{\overnouse}{\ensuremath{\overline{\nouse}} }

\newcommand{\overn}{{\ensuremath{\mathfrak{n}}}}
\newcommand{\y}{{\ensuremath{y}}}
\newcommand{\overcore}{\ensuremath{\overline{\core}} }
\newcommand{\overC}{\ensuremath{\overline{C}} }
\newcommand{\kouse}{\ensuremath{\mathcal{K}} }
\newcommand{\overkouse}{\ensuremath{\overline{\kouse}} }

\renewcommand{\iff}{\textup{\,\ {iff}\,\ }}
\newcommand{\type}{\text{\,--\,type }}
\def\core{{\mathfrak{C}}}
\def\AD{{\textup{AD} }}

\def\DC{{\textup{DC} }}
\def\ZF{{\textup{ZF} }}
\def\PM{{\textup{PM} }}

\def\dom{{\text{dom}}}
\def\ran{{\text{ran}}}

\def\OR{{\textup{OR}}}
\def\BK{{\textup{BK}}}
\def\Hull{{\textup{Hull}}}
\def\Det{{\textup{Det}}}

\newcommand{\mapsigma}[1]{\xrightarrow[\text{ \ \ $\Sigma_{#1}$}]{}}
\newcommand{\maps}[1]{\xrightarrow{\text{#1}}}
 
\def\Rtheory{\textup{R}}
\def\Rplus{\textup{$\Rtheory^+$}}
\def\-{\!-\!}

\newcommand{\boldface}[2]{%
\protect\raisebox{0pt}[0pt][0pt]{%
$\underset{\displaystyle\widetilde{}}{\boldsymbol{#1}}{_{#2}}$}\mbox{\hskip 1pt}}
\newcommand{\boldfaceone}[2]{%
\protect\raisebox{0pt}[0pt][0pt]{%
$\underset{\displaystyle\widetilde{}}{\boldsymbol{#1}}{_{#2}^{1}}$}\mbox{\hskip 1pt}}
\newcommand{\uR}{\underline{\R}}
\newcommand{\uF}{\underline{F}}
\newcommand{\uG}{\underline{\G}}
\newcommand{\up}{\underline{p}}
\newcommand{\urho}{\underline{\rho}}
\newcommand{\uQ}{\underline{Q}}
\newcommand{\ux}{\underline{x}}
\newcommand{\un}{\underline{n}}
\newcommand{\um}{\underline{m}}
\newcommand{\oz}{\hat{z}}
\newcommand{\fn}{\hat{n}}

\newcommand{\cortype}{3.40}
\newcommand{\cortwotype}{3.53}
\newcommand{\deftwo}{3.88}
\newcommand{\Extendible}{3.92}
\newcommand{\goodcovdef}{3.100}
\newcommand{\skolem}{3.72}
\newcommand{\sigmacard}{3.13}
\newcommand{\onemouse}{3.55}
\newcommand{\realpremice}{3.2}
\newcommand{\realmice}{3.4}
\newcommand{\coremice}{3.4.2}
\newcommand{\finestructure}{3}
\newcommand{\minmouseiterable}{3.4.4}
\newcommand{\nicedef}{3.97}
\newcommand{\EOEL}{3.64}
\newcommand{\intismi}{3.89}
\newcommand{\nobigger}{4.4}
\newcommand{\skolemdef}{3.8}
\newcommand{\sounddef}{3.24}
\newcommand{\criterion}{3.93}
\newcommand{\FunnyF}{3.91}
\newcommand{\quasi}{3.98}
\newcommand{\soundlemma}{3.25}
\newcommand{\defontau}{3.99}
\newcommand{\PMdefined}{3.43}
\newcommand{\defontautwo}{3.101}
\newcommand{\soundnplusone}{3.71}
\newcommand{\relationtwo}{3.70}
\newcommand{\mouseitercrit}{3.2.2}

\begin{document}

\title[Scales and the fine structure of $K(\R)$. 
Part II]{Scales and the fine structure of $\boldsymbol{K(\pmb{\R})}$\\ {Part II: Weak real mice and scales}}

\author{Daniel W. Cunningham}
\address{Mathematics Department,
State University of New York,
College at Buffalo,\\
1300 Elmwood Avenue,
Buffalo, NY 14222, USA}
\email{cunnindw@math.buffalostate.edu}

\begin{abstract}
We define weak real mice \ $\mouse$ \ and prove that the
boldface pointclass \ $\boldface{\Sigma}{m}(\mouse)$ \ has the scale property assuming only the
determinacy of sets of reals in \ $\mouse$ \ when \ $m$ \ is the smallest integer \ $m>0$ \ such that
\ $\boldface{\Sigma}{m}(\mouse)$ \ contains a set of reals not in \ $\mouse$. \  We shall use this
development in Part III to obtain scales of minimal complexity in \ $\Kr$.
\end{abstract}

\maketitle

\section{Introduction} 
This paper uses the work presented in Part I \cite{Part1} to address the following question:
\begin{question}[Q]Given an iterable real premouse \ $\mouse$ \ and \ $m\ge 1$, \ when does the boldface\footnotemark
pointclass \ $\boldface{\Sigma}{m}(\mouse)$ \ have the scale property?
\end{question}
Using\footnotetext{We allow arbitrary constants from the domain of the structure $\mouse$.} the fine structure of
real mice presented in \cite{Cfsrm} and
\cite{Part1}, we give a partial answer to this question in section~\ref{finalsection} by proving 
the following theorem  on the existence of scales:
\begin{main}[$\ZF+\DC$]  Suppose that  \ $\mouse$ \  is a
weak real mouse satisfying \ $\AD$. \ Then \
$\boldface{\Sigma}{m}(\mouse)$ \ has the scale property when \ $m=m(\mouse)$. 
\end{main}
The above theorem requires only the determinacy of sets of reals in \ $\mouse$, \ and it extends the following
``lightface'' result established in \cite[see Theorem 4.4]{Crcm}:
\begin{theorem}[$\ZF+\DC$]\label{firstscales}  Suppose that \ $\mouse$ \  is an
iterable real premouse satisfying \ $\AD$. \ Then \ $\Sigma_1(\mouse)$ \
has the scale property. 
\end{theorem}
With Theorem~\ref{firstscales} and Theorem~\ref{newthm} at hand, we will give an explicit answer to
Question (Q) in Part III \cite{Part3}.  

We now give a quick overview leading to the concept of a weak real mouse \ $\mouse$ \ and the definition of
the integer \ $m(\mouse)$. \ We say that \ $\mathcal{M}=(M,\R,\kappa,\mu)$ \ is a real
1--mouse (see \cite[section~{\finestructure}]{Part1}),  if \ $\mouse$
\ is an iterable real premouse and \ $\mathcal{P}(\R\times\kappa)\cap \boldface{\Sigma}{1}(M)\not\subseteq M$, \ where \ $M$ \
has the form \ $J_\alpha[\mu](\R)$ \ and \ $\kappa$ \ is the ``measurable cardinal'' in \ $\mouse$. \ Real 1--mice
suffice to define the real core model and to prove the results in \cite{Crcm} about \ $\Kr$;
\ however, real 1--mice are not sufficient to construct scales of minimal complexity. Our solution
to the problem of identifying these scales in \ $\Kr$ \ requires the development of a {\sl full\/}
fine structure theory for \ $\Kr$. \ In \cite{Cfsrm} we initiated this development by
generalizing the Dodd-Jensen notion of a mouse to that of a {\it real mouse\/}. This is accomplished by (a)
extending the Dodd-Jensen concept of acceptability to include the set of reals, (b) replacing \ $\Sigma_1$ \
with \ $\Sigma_n$, \ where \ $n$ \ is the smallest integer such that \ 
$\mathcal{P}(\R\times\kappa)\cap\boldface{\Sigma}{n+1}(\mouse)\not\subseteq M$, \ and (c) defining an 
iteration procedure stronger than the one we defined in \cite{Crcm}. Now, let \ $\mouse$ \ be a real mouse. Assume that there
is an integer \ $m\ge 1$ \ such that \ $\mathcal{P}(\R)\cap\boldface{\Sigma}{m}(\mouse)\not\subseteq M$ \ and let
\ $m=m(\mouse)$ \ be the least such integer. We say that \ $\mouse$ \ is {\it weak\/}  if 
\be
\item[(1)] $\mouse$ \ is a proper initial segment of an iterable real premouse, and 
\item[(2)] $\mouse$ \ realizes a \ $\Sigma_m$ \ type not realized in any proper initial segment of \ $\mouse$. 
\ee
In (2), a \ $\Sigma_m$ \ type \ $\Upsilon$ \ is a non-empty subset of \[\{ \theta\in \Sigma_m\cup\Pi_m  :
\text{$\theta$ is a formula of one free variable}\}\] and \ $\mouse$ \ is said to {\it realize\/} \ $\Upsilon$
\ if there is an \ $a\in M$ \ such that \ $\mouse\models\theta(a)$ \ for all \ $\theta\in\Upsilon$.

Since the proof of Theorem~\ref{newthm} relies heavily on the fine structure of real mice, 
the proof is more technically involved than the argument used to establish Theorem~\ref{firstscales}. 
So, in sections~\ref{wrm}--\ref{finalsection} (respectively), we shall 
\bi
\item formally define the concept a weak real mouse, 
\item outline the principle idea behind the proof of Theorem~\ref{newthm}, 
\item discuss closed game representations and the construction of scales,
\item investigate the structural properties enjoyed by the core of 
a weak real mouse and define the specific closed games used in the proof of Theorem~\ref{newthm},
\item present a formal proof of Theorem~\ref{newthm}.
\ei

In \cite{Part1} we present our development of the fine structure theory for \ $\Kr$ \ which
will be used in our proof of Theorem~\ref{newthm}. Consequently, we shall presume that the reader 
has access to \cite{Part1}.

\subsection*{Preliminaries and notation}

Let \ $\omega$ \ be the set of all natural numbers.  \ $\R = {^\omega \omega}$ \ is
the set of all functions from \ $\omega$ \ to \ $\omega$. \ We call \ $\R$ \ the set of
reals and regard \ $\R$ \ as a topological space by giving it the product topology,
using the discrete topology on \ $\omega$. \ For a set \  $A \subseteq \R$ \ we
associate a two person infinite game on \ $\omega$, \ with {\it payoff} \
$A$, \ denoted by \ $G_A$:
\[
\begin{aligned}[c]
{}&{\mathbf{I}}\phantom{{\mathbf{I}}} \qquad x(0) \qquad \phantom{x(1)} \qquad x(2) \qquad \phantom{x(3)}\quad \\
{}&{\mathbf{I}}\mathbf{I} \qquad \phantom{x(0)}\qquad x(1) \qquad \phantom{x(2)} \qquad x(3) \quad 
\end{aligned}
\begin{gathered}[c]
{\cdots}
\end{gathered}
\] 
in which player $\mathbf{I}$ wins if \ $x \in A$, \ and \ $\mathbf{II}$ \ wins if \ $x
\notin A$. \ We say that \ $A$ \ is {\it determined\/} if the corresponding game \ $G_A$ \
is determined, that is, either player $\mathbf{I}$ or  $\mathbf{II}$ has a winning
strategy (see \cite[p. 287]{Mosch}). 
The {\it axiom of determinacy\/} \ ($\AD$) \ is a regularity hypothesis about games on \
$\omega$ \ and states: \ $\forall A \subseteq \R \  ( A  \text{\ is
determined})$.

We work in \  $\ZF$ \ and state our additional
hypotheses as we need them.  We do this, in part, to keep a close watch on the use of determinacy in the
proofs of our main theorems. Variables \ $x, y, z, w \dots$ \ generally range over \ $\R$, \ while
\ $\alpha, \beta, \gamma, \delta \dots$ \ (with few exceptions) range over \ $\OR$, \ the class of ordinals.  For \ $x
\in \R$ \ and \ $i\in\omega$ \ we write \ $\lambda.nx(n+i)$ \ for the real \ $y$ \ such that \
$y(n)=x(n+i)$ \ for all  $n$, and we write \ $(x)_i$, \ or \ $x_i$ \ when the context is clear, \ for the real  \ $z$ \ such
that \ $z(n)=x(\seq{n, i })$, \ where \ $\seq{\ , \, }$ \ recursively encodes a
pair of integers by a single integer.  In addition, for \ $x\in\R$ \ and \ $n\in\omega$ \ we write \ $x\restriction
n=\seq{x(0),\dots,x(n-1)}$. \ If \ $0\leq j\leq\omega$ \ and \ $1\leq k\leq\omega$, \ then \ $\omega^j
\times (^\omega\omega)^k$ \ is recursively homeomorphic to \ $\R$, \ and we will implicitly identify the two.  The
cardinal \ $\Theta$ \ is the supremum of the ordinals which are the surjective image of \ $\R$. \ 
For \ $F,G \in [\OR]^{<\omega}$ \ let
\ $F<_{\BK} G \iff \exists \alpha \in G(G=F-\alpha) \lor \max(G \bigtriangleup F)
\in G$. \ Here, \ $\bigtriangleup$ \ is the symmetric difference operation.
The order \ $<_{\BK}$ \ is the Brouwer-Kleene order on finite sets of ordinals and is a \
$\Sigma_0$ \ well-order.

A {\it pointclass\/} is a set of subsets of \ $\R$ \ closed under recursive substitutions. A
boldface pointclass is a pointclass closed under continuous substitutions.  For a
pointclass \ $\Gamma$, \ we write \ ``$\Gamma\-\AD$'' or ``$\Det(\Gamma)$'' \ to denote
the assertion that all games on \ $\omega$ \ with payoff in \ $\Gamma$ \ are
determined. For the concepts of a scale and of the scale property (and any other
notions from Descriptive Set Theory that we have not defined), we refer the reader
to Moschovakis \cite{Mosch}.

A proper class \ $M$ \ is called an {\it inner model\/} if and only if \ $M$ \ is a
transitive  $\in$--model of \ $\ZF$ \ containing all the ordinals. We distinguish between the notations \ $L[A]$ and $L(A)$. \ 
The inner model \ $L(A)$ \ is defined to be the class of sets constructible {\it above\/} \ $A$, \ that is, \ one starts
with a set \ $A$ \ and iterates definability in the language of set theory. \ Thus, \ $L(A)$ \ is  the smallest
inner model \ $M$ \ such that \ $A\in M$. \ The inner model \ $L[A]$ \ is defined to be the class of sets constructible {\it
relative\/} to \ $A$, \ that is, \ one starts with the empty set and iterates definability in the language of set theory
augmented by the predicate \ $A$. \ Consequently, \ $L[A]$ \ is the smallest inner model \ $M$ \ such that \ $A\cap M\in
M$ \ (see page 34 of \cite{Kana}). Furthermore, one defines \ $L[A,B]$ \ to be the class of sets constructible {\it relative\/}
to \ $A$ \ and \ $B$, \ whereas \ $L[A](B)$ \ is defined as the class of sets constructible {\it relative\/} to \ $A$ \ and 
{\it above\/}
\ $B$. \ Thus, \ $A\cap L[A](B)\in L[A](B)$ \ and \ $B\in L[A](B)$.

Our general set theoretic notation is standard. Given a function \ $f$, \ we write
\ $\dom(f) = \{x : \exists y (f(x) =y)\}$ \ and \ $\ran(f) = \{y : \exists x (f(x)
= y)\}$. \ We shall write \ $\seq{x_1,\dots,x_n}$ \ to represent a finite sequence of elements.
For any set \ $X$, \ $(X)^{<\omega}$ \ is the set of all finite sequences of elements
of \ $X$,  \ $[X]^{<\omega}$ \ is the set of all finite subsets of \ $X$, \ and \
$\mathcal{P}(X)$ \ is the set of all subsets of \ $X.$ \ Given two finite sequences \
$s$ \ and \ $t$, \ the sequence \ $s{^\frown}t$ \ is the concatenation of \ $s$ \ to
\ $t$. \ Generally, \ $\mu$ \ will be a normal measure on \ $\mathcal{P}(\kappa)$, \ where
\ $\kappa$ \ is an ordinal. For any ordinals \ $\eta \leq \alpha$, \ 
${^\eta\alpha}\!\uparrow$ \ is the set of all strictly increasing \ $\eta$ \ sequences
from \ $\alpha$. \ $V_\alpha$ \ is the set of all sets of rank less than \ $\alpha$. \ 
We let \ $y = T_c(x)$ \ denote the formula ``$y$ \ is the transitive closure of \ $x$.'' \ 
 For a model \ $\mouse = (M,\in, \dots),$ \ we shall abuse standard notation slightly and write
 \ ${^\kappa M}=\{f\in M \ \vert \ f\colon \kappa \rightarrow M \}$. \ In addition, for a model (or inner model) \
$\mouse$ \ having only one ``measurable cardinal,'' we shall write \ $\kappa^{\mouse}$ \ to denote this cardinal
in \ $\mouse$. \ Similarly, when \ $\mouse$ \ has only one ``measure,'' we shall write \ $\mu^{\mouse}$ \ to denote
this measure.

Given a model \ $\mouse = (M,c_1,c_2,\dots,c_m,A_1,A_2,\dots,A_N),$ \ where the \  $A_i$ \
are predicates and the \ $c_i$ \ are constants,  if \ $X\subseteq M$ \ then \
$\Sigma_n(\mouse,X)$ \ is the class of relations on \ $M$ \ definable over \ $\mouse$ \ by a
\ $\Sigma_n$ \ formula from parameters in \ $X \cup \lbrace c_1,c_2,\dots,c_m\rbrace$. \
$\Sigma_\omega(\mouse,X) =  \bigcup\limits_{n\in\omega} \Sigma_n(\mouse,X).$  We write
``$\Sigma_n(\mouse)$'' for \ $\Sigma_n(\mouse,\emptyset)$ \ and \
``$\boldface{\Sigma}{n}(\mouse)$'' for  the boldface class \ $\Sigma_n(\mouse,M).$ \ 
Similar conventions hold for \ $\Pi_n$ \ and \ $\Delta_n$ \ notations.  If \ $\mouse$ \ is
a substructure of \ $\nouse$ \ and \ $X\subseteq M \subseteq N$, \ then 
``$\mouse \prec_n^X \nouse$'' means that \ $\mouse \models \phi[a] \text{ \ iff
 \ } \nouse \models \phi[a]$, \ for all \ $a \in (X)^{<\omega}$ \ and for all \
$\Sigma_n$ \ formulae \ $\phi$ \ (the formula \ $\phi$ \ is allowed constants taken
from \ $\lbrace c_1,c_2,\dots,c_m\rbrace$). \ We write \ ``$\mouse \prec_n \nouse$'' \ 
for \ ``$\mouse \prec_n^M \nouse$.'' \ Also, for any two models \ $\mouse$ \ and \
$\nouse$, \  we write \ $\pi:\mouse\mapsigma{n}\nouse$ \ to indicate that the map \ $\pi$ \ is
a \ $\Sigma_n$--elementary embedding, that is, \ $\mouse \models \phi[a] \text{ \ iff
 \ } \nouse \models \phi[\pi(a)]$, \ for all \ $a=\seq{a_0,a_1,\dots} \in (M)^{<\omega}$ \
and for all \ $\Sigma_n$ \ formulae \ $\phi$, \ where \ $0\le n\le\omega$ \ and \ 
$\pi(a)=\seq{\pi(a_0),\pi(a_1),\dots}$.

\section{Weak real mice}\label{wrm}
In \cite[subsection~{\realpremice}]{Part1}  we defined the notion of a
premouse `above the reals' in the language \ $\Lng_n$, \ where \ $n\ge 0$. \ In \cite[subsection
{\realmice}]{Part1} we defined the concept of a mouse `above the reals'. Our objective now is to define a
{\it weak mouse above the reals\/} (see Definition~\ref{weakmouse}).
\begin{definition}Let \ $\mouse$ \  be a mouse and let \ $n=n(\mouse)$, \ that is, \ $\rho_{\mouse}^{n+1}\le
\kappa^\mouse < \rho_{\mouse}^n$. \ Define \ $m=m(\mouse)$ \ to be the least integer \ $m\ge n+1$
\ such that \ $\rho_{\mouse}^m=1$, \ if such an \ $m$ \ exists. Recall that for such an \ $m$, \ ${\mathcal
P}(\R)\cap\boldface{\Sigma}{m}(\mouse)\not\subseteq M$. 
\end{definition} 
\begin{definition}\label{type} Suppose \ $\mouse$ \ is a premouse in the language \ $\Lng_N$
\ and let \ $a\in M$. \ The \ $\Upsilon_n^N$--{\it type realized by \ $a$ \ in \ \mouse} \ is the set
\[\Upsilon_n^N(a,\mouse) = \{\,\vartheta(v) \in \Sigma_n^N \cup \Pi_n^N:  \mouse\models
\vartheta[a]\,\}\]
 where \ $\vartheta(v)$ \  is an \ $\Lng_N$ \ formula with one free
variable. Let \ $\Upsilon \subset \Sigma_n^N \cup \Pi_n^N$ \ be nonempty and consist of formulae with one free
variable. The type \ $\Upsilon$ \ is {\it realized in\/} \ $\mouse$ \ if for
some \ $b\in M$, \ $\Upsilon \subseteq \Upsilon_n^N(b,\mouse)$.
\end{definition}
We now define the notion of a weak mouse. 
\begin{definition}\label{weakmouse} Let \ $\mouse$ \  be a mouse. \ Suppose that \ $m=m(\mouse)$
\ is defined and that \ $\mouse$ \  is a proper initial segment of an iterable pure premouse. 
Then  \ $\mouse$ \ is said to be a {\it weak mouse} if there is a \ $\Upsilon_m^0$ \ type realized
in  \ $\mouse$ \  which is not realized in any proper initial segment of \ $\mouse$, \ that is, for some \ $a\in M$
\[\forall \gamma < \widehat{\OR}^\mouse, \ \mouse^\gamma \text{ does not realize }  \Upsilon_m^0(a,\mouse).\]
In contrast, if every \
$\Upsilon_m^0$ \ type realized in  \ $\mouse$ \  is also realized in a proper initial
segment of \ $\mouse$, \ then \ $\mouse$ \ is said to be a {\it strong mouse.}
 \end{definition}
The following lemma is used in section~\ref{cwrm} to identity a connection between the types realized by a weak real
mouse and the types realized by its core.
\begin{lemma}\label{lowrho} Let \ $\mouse$ \ be a weak mouse with \ $n=n(\mouse)$. \ Then \ $\rho_{\mouse}^{n+1}<
\kappa^\mouse$.
\end{lemma}
\begin{proof} Since \ $\mouse$ \ is weak, we have that \ $m=m(\mouse)$ \ is defined and  \ $\rho_{\mouse}^m=1$.
\ Corollary~{\nobigger} of \cite{Part1} now implies that \ $\rho_{\mouse}^{n+1}<\kappa^\mouse$.
\end{proof}

\section{\bf The fundamental idea behind the proof of Theorem~\ref{newthm}}\label{sectwo}
We now give some motivation behind the proof Theorem~\ref{newthm}. First we recall the main idea
supporting the proof of  Theorem~\ref{firstscales} in \cite[see pages 268-278]{Crcm}. Theorem~\ref{firstscales} asserts
that if \ $\mouse$ \  is an iterable real premouse satisfying \ $\AD$, \ then any \ $\Sigma_1(\mouse)$ \ set of reals
\ $P$ \ has a scale which is also \ $\Sigma_1(\mouse)$. \ Let \ $\mouse$ \  be an iterable premouse and let \
$P$ \ be \ $\Sigma_1(\mouse)$ \ set of reals.  Define   \[F^{\mouse}=\{\,f\in M :  \exists n\in\omega \ \mouse
\models f\colon {^n}\underline\kappa \rightarrow \OR\,\}.\] For \ $f\in F^{\mouse}$, \ write \ $d(f)
= n$ \ if and only if \ $n\in\omega$  \ and \ $\mouse \models f\colon {^n}\underline\kappa
\rightarrow \OR$. \ We shall assume the convention that \ $f\in F^{\mouse}$ \ and  $d(f) = 0$ \
whenever \ $f\in {\OR}^{\mouse}$. \ Finally, for \ $n\in\omega$, \ define
\  $F^{\mouse}_n = \{\,f\in F^{\mouse} : d(f) = n \,\}$. \  We shall now review the strategy behind our proof
of Theorem~\ref{firstscales} in \cite{Crcm}. The key idea in the proof was to design a closed game
representation \ $x\mapsto G_x$ \ for \ $P$ \ as follows: For each \ $x\in\R$ \ there is a game \ $G_x$ \ in
which player $\mathbf{I}$'s  moves come from \
$\R\times F^{\mouse}$ \ while player \ $\mathbf{II}$'s \ moves come from
\ $\R$. \ Thus, a typical run of the game \ $G_x$ \ has the form 
\[
\begin{aligned}[c]
{}&{\mathbf{I}}\phantom{{\mathbf{I}}} \qquad x_0,f_0 \qquad \phantom{x_1} \qquad x_2,f_1 \qquad
\phantom{x_3}\quad \\ {}&{\mathbf{I}}\mathbf{I} \qquad \phantom{x_0,f_0}\qquad x_1 \qquad \phantom{x_2,f_1}
\qquad x_3 \quad 
\end{aligned}
\begin{gathered}[c]
{\cdots}
\end{gathered}
\] 
where \ $x_i\in\R$ \ and \ $f_i\in F_{c(i)}^{\mouse}$, \ for some recursive function\
$c\colon\omega \rightarrow \omega$.
\ The game \ $G_x$ \ is {\it closed and continuously associated to} \ $x$ \ if
for some \ $Q\subseteq (\omega^{<\omega})^{<\omega}\times(F^{\mouse})^{<\omega}$, \ the
following hold: 
\be
\item For each \ $n\in\omega$ \ the relation \ $Q_n$ \ is invariant.
\item Player $\mathbf{I}$ wins \ $G_x$ \ if and only if \ $\forall n
Q_n(x,x_0,\dots,x_n,f_0,\dots,f_n)$, \ where for each \ $n\in\omega$ 
\[Q_n = \{\,(x,x_0,\dots,x_n,f_0,\dots,f_n) : 
Q(\seq{x\restriction n,x_0\restriction n,\dots,x_n\restriction n},\seq{f_0,\dots,f_n})\,\}.\]
\ee

By Gale-Stewart (see \cite[pages 289, 446-7]{Mosch}), one of the players has a winning
quasi-strategy. Since we are not assuming the axiom of choice, we do not get single-valued
strategies. 

\begin{definition}\label{oldgames} For a set \ $P\subseteq\R$, \ we say that \ {\it $P$ \ admits
a closed game representation,} if there is a map \ $x \mapsto G_x$ \ such that, for all \ $x\in\R$
\[P(x) \Longleftrightarrow \mathbf{I} \text{ wins } G_x,\]
where \ $G_x$ \ is closed and continuously associated to \ $x$.
\end{definition}
Suppose that \ $P\subseteq\R$ \ admits a closed game representation \ $x \mapsto G_x$. \  Let
\begin{alignat*}{1}
P_k(x,u) \iff &u \text{ \ is a position in \ }G_x \text{ \ of length \ }k \text{ \ from
which } \tag{$*$}\\ 
&\text{player }\mathbf{I}\text{ \ has a winning quasi-strategy.}
\end{alignat*}
Here, \ {\it $u$ \ is a position in \ $G_x$ \ of length \ $k$} \ if \ $u$ \ has 
the form \ $u = \seq{\seq{x_{2i},f_i,x_{2i+1}} \ : 0\le i < k}$. \ We can then build a
scale on \ $P$ \ using a generalization of a scale construction due to Moschovakis
\cite{Moone}.  One first defines a scale on \ $P_k$ \ for all \ $k$ \ simultaneously,
assuming the axiom of dependent choice and enough determinacy. We are then able to construct a scale on \
$P$, \ using the scales on each \ $P_k$. 

Let \ $\mouse$ \ be an iterable real premouse and let \ $P$ \ be a set of reals definable by a \ $\Sigma_1$ \ formula  \
$\varphi(v)$ \ in a proper initial segment of \ $\mouse$. \ In our proof of Theorem~\ref{firstscales} in \cite{Crcm}, we 
designed a specific closed game representation \ $x\mapsto G_x$ \  for \ $P$ \ (see \cite[p. 266]{Crcm}) simple enough to
ensure that 
\bi
\item the associated \ $P_k\in\mouse$ \ for all \ $k\in\omega$,
\item the Moschovakis scale on \ $P$ \ is \ $\Sigma_1(\mouse)$ \ and its construction requires only the determinacy of sets
of reals in \ $\mouse$.
\ei
The basic plan  behind the design of the game \ $G_x$ \ is to force player $\mathbf{I}$ (if he wants to win) to
construct an iterable model which contains all the reals played in the run of the game.  The game's payoff is
defined so that the winning player $\mathbf{I}$ must construct a premouse \ $\nouse$ \ in which \ $\varphi(x)$ \ holds
and, in addition, player $\mathbf{I}$ must play functions \ $f\in F^{\mouse}$ \ which verify that \ $\nouse$ \ is
{\it premouse iterable\/}. It turns out that the only way for player $\mathbf{I}$ to win this game is to use \
$\Hull_1^{\mouse}(\R)$ \ as a guide in the construction of his iterable model \ $\nouse$ \ (see \cite[Lemma
4.7]{Crcm}). \ Consequently, one obtains the desired closed game representation for \ $P$. \ We make the
following observation:  {\sl The canonical model \ $\Hull_1^{\mouse}(\R)$ \ is the key
ingredient in the construction of the closed game representation for such a set as \ $P$.}

Suppose now that  \ $\mouse$ \  is a  weak real mouse satisfying \ $\AD$. \  Theorem~\ref{newthm} asserts that 
\  $\boldface{\Sigma}{m}(\mouse)$ \ has the scale property where \ $m=m(\mouse)$. \ The ideas supporting the proof of
this theorem is simply stated as follows:  Let \ $P$ \ be a \ $\boldface{\Sigma}{m}(\mouse)$ \ set of reals. To
construct  a \ $\boldface{\Sigma}{m}(\mouse)$ \ scale on \ $P$, \ we shall show
that \ $P=\bigcup\limits_{i\in\omega}P^i$ \ where each \ $P^i\subseteq\R$ \ admits a closed game representation \ $x
\mapsto G^i_x$. \ These closed game representations \ $x\mapsto G^i_x$ \  for \ $P^i$ \ 
are simple enough so that the scale constructed on \ $P$  \  (see section~\ref{cgr} below) is \
$\boldface{\Sigma}{m}(\mouse)$. \ As in the proof of Theorem~\ref{firstscales} discussed above, we shall identify a
canonical model that will allow us to construct the desired closed game representations on each \ $P^i$. \ The
obvious candidate for such a model is \ $\core=\core(\mouse)$, \ the core of \
$\mouse$. \ However, we want a closed game representation in which player $\mathbf{I}$ can  easily verify that the model he
constructs is {\it mouse iterable\/}.  Since the definition of mouse iterability involves the structure \
$\overcore$, \ we shall use \ $\overcore$ \ as the canonical model in our definition of the closed game representation
for each \ $P^i$ \ and we will require a winning player $\mathbf{I}$ to play elements \ $\f\in\F^{\overcore}$ \ to verify
that his model is mouse iterable (see \cite[section~{\mouseitercrit}]{Part1}). 
\section{Closed game representations and scales}\label{cgr}
Our proof of Theorem~\ref{newthm} constructs the desired scale by means of a closed game
representation. In this section, we shall first describe the kind of closed game that will be
used in our proof of Theorem~\ref{newthm}. We shall then discuss two relevant methods for
constructing scales from such closed games. 

\subsection{Closed game representations}
The key concept  behind the proof of Theorem~\ref{firstscales} in \cite{Crcm} is the notion of a 
``closed game representation'' as described in Definition~\ref{oldgames} above. We need to  modify this concept slightly
for our proof of  Theorem~\ref{newthm}; but before we do this, we introduce some notation. 

\begin{definition} Let \ $\mouse$ \ be a real mouse. We say that \ $\mathcal{X} =
X_1\times\dots\times X_k$ \ is an \ {\it $\overline{\mouse}$--space}, \ if for all \ $i\le k$,
\ either \ $X_i = \R$ \  or \ $X_i = \F^{\overline{\mouse}}_m$ \ where \ $m\in\omega$. \ For \
$u =
\seq{u_1,\dots,u_k}$ \ and \ $u^\prime = \seq{u_1^\prime,\dots,u_k^\prime}$ \ in \
$\mathcal{X}$ \ we write \ $u\sim u^\prime$ \ if for all \ $i\le k$,
\[\begin{cases} u_i = u^\prime_i, &\text{if \ $X_i = \R$} \\
         u_i \equiv_{\mu_m} u^\prime_i, &\text{if \ $X_i = \F^{\overline{\mouse}}_m$ \ for some 
\ $m\in\omega$}.\end{cases}\]
\end{definition}
The set \ $\F^{\overline{\mouse}}_m$ \ in the above definition is described in
\cite[Definition~{\FunnyF}]{Part1}.
\begin{definition} Let \ $\mouse$ \ be a real mouse and let \ $\mathcal{X}$ \ be an \
$\overline{\mouse}$--space. A relation \ $P\subseteq \mathcal{X}$ \ is said to be  {\it invariant\/}  if
\[u\sim u^\prime \implies [P(u) \Leftrightarrow P(u^\prime)]\]
for all \ $u, u^\prime \in \mathcal{X}$.
\end{definition}

\begin{definition}\label{newclsdgame}
Let \ $\mouse$  \ be a real mouse and let \ $c\colon\omega \rightarrow \omega$ \ be a recursive
function. Suppose that for each \ $x\in\R$ \ there is a game \ $G_x$ \ in which player  $\mathbf{I}$'s  moves come from \
$\R\times \F^{\overmouse}$ \ while player \ $\mathbf{II}$'s \ moves come from
\ $\R$. \ Thus, a typical run of the game \ $G_x$ \ has the form 
\[
\begin{aligned}[c]
{}&{\mathbf{I}}\phantom{{\mathbf{I}}} \qquad x_0,\f_0 \qquad \phantom{x_1} \qquad x_2,\f_1 \qquad \phantom{x_3}\quad \\
{}&{\mathbf{I}}\mathbf{I} \qquad \phantom{x_0,\f_0}\qquad x_1 \qquad \phantom{x_2,\f_1} \qquad x_3 \quad 
\end{aligned}
\begin{gathered}[c]
{\cdots}
\end{gathered}
\] 
where \ $x_i\in\R$ \ and \ $\f_i\in
\F_{c(i)}^{\overmouse}$.
\ We shall say that the game \ $G_x$ \ is {\it closed and continuously associated to} \ $x$ \ if
for some \ $Q\subseteq (\omega^{<\omega})^{<\omega}\times(\F^{\overmouse})^{<\omega}$, \ the
following hold: 
\be
\item For each \ $n\in\omega$, \ the relation \ $Q_n$ \ is invariant.
\item Player $\mathbf{I}$  wins \ $G_x$ \ if and only if \ $(\forall n)
Q_n(x,x_0,\dots,x_n,\f_0,\dots,\f_n)$, \ where for each \ $n\in\omega$
\[Q_n = \{\,(x,x_0,\dots,x_n,\f_0,\dots,\f_n) : 
Q(\seq{x\restriction n,x_0\restriction n,\dots,x_n\restriction n},\seq{\f_0,\dots,\f_n})\,\}.\]
\ee
\end{definition}
Hence by Gale-Stewart (see \cite[pp. 289, 446-7]{Mosch}), one of the players has a winning
quasi-strategy. Since we are not assuming the axiom of choice, we do not get single-valued
strategies. 

\begin{definition} Let \ $\mouse$ \ be a real mouse. For a set \ $P\subseteq\R$, \ we say that \ {\it $P$ \ admits
a closed game representation via $\mouse$,} \ if there is a map \ $x \mapsto G_x$ \ such that, for all \ $x\in\R$
\[P(x) \iff \mathbf{I} \text{ wins } G_x,\]
where \ $G_x$ \ is closed and continuously associated to \ $x$ \ as in Definition~\ref{newclsdgame}.
\end{definition}

\subsection{The Moschovakis scale}
Suppose that \ $P\subseteq\R$ \ admits a closed game representation \ $x \mapsto G_x$ \ via the real mouse \ $\mouse$. \ We
can build a scale on \ $P$ \ using a generalization of a scale construction due to Moschovakis
\cite{Moone}.  Let
\begin{gather}\label{partial_play}
\begin{split}
P_k(x,u) \iff &\, u \text{ is a position in } G_x \text{ of length } k \text{ from
which player $\mathbf{I}$}\\ 
&\text{ \ has a winning quasi-strategy.}
\end{split}
\end{gather}
Here, \ {\it $u$ \ is a position in \ $G_x$ \ of length \ $k$} \ if \ $u$ \ has 
the form \ $u = \seq{\seq{x_{2i},\f_i,x_{2i+1}} \ : 0\le i < k};$ \ 
however, we shall abuse this notation slightly and identify \ $u$ \ with the sequence
\ $u=\seq{ x_{0},\f_0,x_{1},x_{2},\f_1,x_{3},\dots,x_{2(k-1)},\f_{k-1},x_{2(k-1)+1}}$.

Note that \ $P(x) \Leftrightarrow P_0(x,\emptyset)$ \ and for each \ $k\in\omega$, \ $P_k$
\ is an invariant relation. We extend the concept of scale to the relations \ $P_k$ \ by
giving \ $\F_n^{\mouse}$, \ for each \ $n\in\omega$, \ the \ $\equiv_{\mu_n}$--discrete topology,
that is, \ $\f_i \rightarrow \f$ \ if \ $\exists m\forall i\ge m \ (\f_i \equiv_{\mu_n} \f).$ \ Following
Moschovakis we define scales \ $\seq{\varphi_i^k : i \in \omega}$ \
on \ $P_k$ \ for all \ $k$ \ simultaneously, assuming the axiom of dependent choice and enough
determinacy. 
First, note that by (\ref{partial_play}) above we have that
\begin{equation*}P_k(x,u) \iff (\exists w \in \R)(\exists \f\in \F^{\overmouse}_{c(k)})(\forall y\in \R)
P_{k+1}(x,u{^\frown}\seq{w,\f,y}).\end{equation*}
Define the intermediate relations:
\begin{align*}
\Hat{P}_k(x,u;w) \iff &(\exists \f\in \F^{\overmouse}_{c(k)})(\forall y\in \R)
P_{k+1}(x,u{^\frown}\seq{w,\f,y})\\
\Hat{\Hat{P}}_k(x,u;w,\f) \iff &(\forall y\in \R)
P_{k+1}(x,u{^\frown}\seq{w,\f,y}).
\end{align*}
Notice that 
\begin{align*}
P_k(x,u) \iff &\, (\exists w \in \R)\,\Hat{P}_k(x,u;w)\\
\Hat{P}_k(x,u;w) \iff &\, (\exists \f\in \F^{\overmouse}_{c(k)})\,\Hat{\Hat{P}}_k(x,u;w,\f).
\end{align*}

We define scales \ $\seq{\varphi_i^k : i \in \omega}$, 
$\seq{\Hat{\varphi}_i^k : i \in \omega}$, 
$\seq{\Hat{\Hat{\varphi}}_i^k : i \in \omega}$ \ 
on \
$P_k,\,\Hat{P}_k,\,\Hat{\Hat{P}}_k$, \ respectively, by induction. 
If \ $P_k(x,u)$, \ then \
$\varphi_0^k(x,u) = 0$; \ otherwise \ 
$\varphi_0^k(x,u)$ \ is undefined. Define \ $\Hat{\varphi}_0^k$ \ and \ $\Hat{\Hat{\varphi}}_0^k$ \
similarly. Now define \ 
\be
\item $\varphi_{i+1}^k$ \ from \ $\Hat{\varphi}_0^k,\dots, \Hat{\varphi}_i^k$
\item $\Hat{\varphi}_{i+1}^k$ \ from \ $\Hat{\Hat{\varphi}}_i^k$
\item $\Hat{\Hat{\varphi}}_{i+1}^k$ \ from \ $\varphi_0^{k+1},\dots, \varphi_i^{k+1}$
\ee
by considering the possible cases, respectively.

{\bf Case 1:} $P_k(x,u)$. \ Then \ $(\exists w \in \R)\Hat{P}_k(x,u;w)$ \ and define \ $\varphi_{i+1}^k$ \ to be 
``inf''$\{\Hat{\varphi}_0^k, \dots , \Hat{\varphi}_i^k\}$, \ as in 
\cite{Moone},
so that
\[
\varphi_{i+1}^k(x,u) =\textup{infimum}\{ \seq{\Hat{\varphi}_0^k(x,u;w), w(0),\dots ,\Hat{\varphi}_i^k(x,u;w), w(i)} :
\Hat{P}_k(x,u;w)\}.
\]
Here, \ $\seq{\Hat{\varphi}_0^k(x,u;w), w(0), \dots ,\Hat{\varphi}_i^k(x,u;w), w(i)}$ \ is
the ordinal of this tuple in the lexicographic order.

{\bf Case 2:} $\Hat{P}_k(x,u;w)$. \ Then \ $(\exists
\f\in \F^{\overmouse}_n)\Hat{\Hat{P}}_k(x,u;w,\f)$ \ letting \ $n={c(k)}$. \ Define  \ $\Hat{\varphi}_{i+1}^k$ \ to be
``min''$\{\Hat{\Hat{\varphi}}_i^k\}$ \ so that
\[\Hat{\varphi}_{i+1}^k(x,u;w) 
= \seq{|\, [\f]_{\mu_n}|,\Hat{\Hat{\varphi}}_i^k(x,u;w,\f)}\]
where \ $|\ [\f]_{\mu_n}|$ \ denotes the \
$E_{\mu_n}^{\overmouse}$--rank \ of \ $\f$ \ (see \cite[Definition~{\nicedef}]{Part1}) \ and \ $\f$ \ is such that
\be
\item $\Hat{\Hat{P}}_k(x,u;w,\f)$
\item $(\forall \g \in \F_n^{\overmouse})(\,\g \,E_{\mu_n}\, \f \Rightarrow \neg
\Hat{\Hat{P}}_k(x,u;w,\g)\,)$.
\ee

{\bf Case 3:} $\Hat{\Hat{P}}_k(x,u;w,\f)$. \ Then \ $(\forall y\in\R)P_{k+1}
(x,u{^\frown}\seq{w,\f,y})$. \ Define \ $\Hat{\Hat{\varphi}}_{i+1}^k$ \ to be
``fake sup''$\{\varphi_0^{k+1},\dots,\varphi_i^{k+1}\}$, \ that is, the ``fake supremum norm''
associated with \ $\varphi_0^{k+1},\dots,\varphi_i^{k+1}$. \ This norm
is defined in detail in \cite{Moone} and its construction uses determinacy of the so-called {\it sup
games.\/} For the benefit of the reader who may not be familiar with \cite{Moone}, we give a
brief overview of the construction of this fake supremum norm. Let \
$\tau_0,\tau_1,\dots$ \ be a fixed enumeration all finite sequences of \ $\omega$, \ so that \
$\tau_0 =\emptyset$ \ and if \ $\tau_m$ \ is a proper initial segment of \ $\tau_j$, \ then \
$m<j$. \ For
\ $(x,u;w,\f),\, (x^\prime,u^\prime;w^\prime,\f^\prime) \in \Hat{\Hat{P}}_k$ \ define \
$(x,u;w,\f)\le^*(x^\prime,u^\prime;w^\prime,\f^\prime)$ \ if and only if \ player $\mathbf{II}$ has a
winning strategy in the following game  \ $G=G(x,u,w,\f,x^\prime,u^\prime,w^\prime,\f^\prime)$ \ on \ $\omega$:
\begin{equation*}
\begin{aligned}[c]
{}&{\mathbf{I}}\phantom{{\mathbf{I}}} \qquad z(0) \qquad \phantom{z'(0)} \qquad z(1) \qquad \phantom{z'(1)}\quad \\
{}&{\mathbf{I}}\mathbf{I} \qquad \phantom{z(0)}\qquad z'(0) \qquad \phantom{z(1)} \qquad z'(1) \quad 
\end{aligned}
\begin{gathered}[c]
{\cdots}
\end{gathered}
\end{equation*}
and, letting \ $y=\tau_i{^\frown}z$ \ and \ $y^\prime=\tau_i{^\frown}z^\prime$, \ player $\mathbf{II}$ wins \ if
and only if
\begin{equation}\begin{aligned}
\langle \varphi_0^{k+1}(x,u{^\frown}\langle w,\f,y\rangle),
 \dots ,&\varphi_i^{k+1}(x,u{^\frown}\langle w,\f,y\rangle)\rangle 
\le \\
\langle \varphi_0^{k+1}&(x^\prime,u^\prime{^\frown}\langle w^\prime,\f^\prime,y^\prime\rangle),
 \dots ,\varphi_i^{k+1}(x^\prime,u^\prime{^\frown}\langle
w^\prime,\f^\prime,y^\prime\rangle)\rangle
\end{aligned}\tag{$\blacktriangle$}\end{equation}
where \ $\le$ \ is the lexicographic order on tuples of ordinals. Let \
$A=A(x,u,w,\f,x^\prime,u^\prime,w^\prime,\f^\prime)$  \ denote the payoff set, defined by the above $(\blacktriangle)$, for the game \ $G$. \ 
Assuming DC and the determinacy of each set of reals \ $A$, \ one can show that the relation \
$\le^*$ \ is a prewellordering on \ $\Hat{\Hat{P}}_k$. \ Define \ $\Hat{\Hat{\varphi}}_{i+1}^k(x,u;w,\f)$ \ to be the \
$\le^*$--rank of \
$(x,u;w,\f)$.

We remark that 
\[[P_k(x,u) \land u\sim u^\prime] \implies \varphi_i^k(x,u) = \varphi_i^k(x,u^\prime).\]
As in \cite{Moone}, assuming enough determinacy, one can show that each \ 
$\seq{\varphi_i^k : i \in \omega}$ \ is a scale on \ $P_k$ \ and thus, \ 
$\seq{\varphi_i^0 : i \in \omega}$ \ is a scale on \ $P$ \ as desired.

\begin{lemma}[$\AD+\DC$] Suppose that \ $P\subseteq\R$ \ admits a closed game representation via a real
mouse \ $\mouse$. \ Then \ $\seq{\varphi_i^0 : i \in \omega}$ \ is a scale on \ $P$.
\end{lemma}

\begin{definition} For \ $P$ \ as above, we call \ $\seq{\varphi_i^0 : i \in \omega}$ \ the {\it Moschovakis scale} \ on \
$P$.
\end{definition}

The prewellordering \ $\le_i$, \ induced by the \ $i^{\underline{\text{th}}}$ \ norm \
$\varphi_i^0$ \ in this scale on \ $P$, \ is constructed from the \ $P_k$'s \ for \ $k\le i$ \ by
means of rudimentary operations, including quantification over \ $\R$ \ and quantification over \
$\F_n^{\overmouse}$ \ for finitely many \ $n\in\omega$. \ Thus, the amount of determinacy required
to construct the scale is closely related to the definability of the scale constructed.
Furthermore, if one can show that each set of reals \ $A(x,u,w,\f,x^\prime,u^\prime,w^\prime,\f^\prime)$ \ defined by
$(\blacktriangle)$ is in \ $\overmouse$, \ then \ $\Det(\overmouse \cap \pow(\R))$ \ is sufficient to conclude that \ $\seq{\le_i\, :
i\in\omega}$ \ is a scale on \ $P$. 
The following lemma will be implicitly applied in this paper. 
\begin{lemma}[$\ZF+\DC$]\label{poffinmouse} Let  \ $\mouse$ \ be a real mouse such that all sets of reals in \ $\overmouse$ \
are determined. Suppose that \ $P\subseteq\R$ \ admits a closed game representation via \ $\mouse$ \ such that each set
of reals as defined by $(\blacktriangle)$ is in \ $\overmouse$. \ Then
\ $\seq{\varphi_i^0 : i \in \omega}$ \ is a scale on \ $P$.
\end{lemma}
\subsection{The union scale}\label{unionscale}
Suppose that \ $P\subseteq\R$ \ is such that \ $P=\bigcup\limits_{i\in\omega}P^i$ \ where each 
\ $P^i\subseteq\R$ \ admits a closed game representation \ $x \mapsto G^i_x$. \ We shall build a
scale on \ $P$. \ Since each \ $P^i\subseteq\R$ \ admits a closed game representation, we shall
let \ $\seq{\varphi_j^i : j \in \omega}$ \ be the {\it Moschovakis scale} \ on \ $P^i$. \ We are assuming
that the relevant games in the construction of these scales are determined, and so we can now define a
scale \ 
$\seq{\psi_j : j \in \omega}$ \ on \ $P$ \ as follows: 
\begin{equation}\label{uscale}
\begin{aligned}
\psi_0(x) &= \text{the least \ $i\in\omega$ \ such that \ $P^i(x)$}\\
\psi_{j+1}(x) &= \seq{\psi_0(x), \varphi_j^{\psi_0(x)}(x)}
\end{aligned}
\end{equation}
where \ $\psi_{j+1}$ \ is defined using the lexicographic order to assign ordinals to pairs of ordinals.
It is not hard to verify that \ $\seq{\psi_j : j \in \omega}$ \ is a scale on \  $P$.
\ The following lemma summarizes these observations.
\begin{lemma}[$\AD+\DC$] Suppose that \ $P\subseteq\R$ \ is such that \
$P=\bigcup\limits_{i\in\omega}P^i$ \ and each \ $P^i\subseteq\R$ \ admits a closed game representation \ $x \mapsto G^i_x$. \ Then \ $\seq{\psi_j : j \in
\omega}$, \ as defined in (\ref{uscale}), is a scale on \ $P$.
\end{lemma} 
\begin{definition}\label{defunionscale} Let \ $P$, \ $\seq{P^i : i\in\omega}$ \ and \ $\seq{\psi_j : j \in
\omega}$ \ be as in the above lemma. We shall call \ $\seq{\psi_j : j \in
\omega}$ \ the {\it union scale\/} on \ $P$.
\end{definition}
\section{The core of a weak mouse}\label{cwrm}
Let \ $\mouse$ \ be a weak real mouse with core \ $\core=\core(\mouse)$. \  We will show in this section that the
core mouse \ $\core$ \ enjoys important structural properties not possessed by a typical core mouse. These additional
properties will allow us to construct our desired closed game representations which will be used to show that  \
$\boldface{\Sigma}{m}(\mouse)$
\ has the scale property, when  \ $m=m(\mouse)$.

Recall that a real mouse \ $\mouse$ \ contains all the reals, that is, \
$\R^\mouse=\R$. \ However, the relevant structural properties that we shall present in this
section hold for weak mice \ $\mouse$ \ in general, and so, until further notice we will only assume that \
$\R^\mouse\subseteq\R$.\footnote{We shall assume that
$\R^\mouse\subseteq\R$ until Remark
\ref{allreals}. The results established prior to this remark shall be applied, as part of the proof of Theorem
\ref{goodcovering}, in a generic extension that adds reals.} 

Throughout this section we let \ $\mouse=(M,{\R}^{\mouse},\kappa,\mu)$ \ be a fixed weak
mouse.  Also fix \ $\core=\core(\mouse)$, \ $\overn=n(\mouse)$ \ and \ $m=m(\mouse)$. \ 
Recall that \ $\overmouse=\mouse^\overn$ \ and \ $\overcore=\core^\overn$. \ Let \ $k\in\omega$ \ be such that \
$m=\overn+k$ \ and note that \ $k\ge 1$. \ Let \ $a\in M$ \ be such that the type \
$\Upsilon=\Upsilon^0_m(a,\mouse)$ \ witnesses that \ $\mouse$ \ is weak. The type \ $\Upsilon$ \ will also be fixed
in this section.

\begin{lemma}\label{T3.10} Assume that \ $\fnouse$ \ is a mouse with \ $n(\fnouse)=\overn$. \ Suppose that \
$\alpha\in\OR$ \ is a multiple of \ $\omega^\omega$. \ Let \ $\fnouse_\alpha$ \ be the \ $\alpha^{\un{\text{th}}}$ \
mouse iterate of \ $\fnouse$. \  Then \ $\fnouse_\alpha$ \ realizes \ $\Upsilon$ \ if and only if \
$\overline{\fnouse}$ \ realizes \ $\overline{\Upsilon}^*$.
\end{lemma}
\begin{proof}
Corollary~{\cortype} of \cite{Part1} implies that \ $\Upsilon$ \ can be translated to the \
$\Upsilon^{\overn}_{k}$\type \ $\overline{\Upsilon}$ \ such that \ $\nouse_\alpha$ \ realizes \ $\Upsilon$
\ if and only if \ $\overnouse_\alpha$ \ realizes \ $\overline{\Upsilon}$. \ Corollary~{\cortwotype} of \cite{Part1}
and Lemma~{\deftwo}  of \cite{Part1}, imply that the type \ $\overline{\Upsilon}$ \ can be translated to the \
$\Upsilon^{\overn}_k$\type \ $\overline{\Upsilon}^*$ \ such that \ $\overnouse_\alpha$ \ realizes \
$\overline{\Upsilon}$
\  if and only if \ $\overnouse$ \ realizes \ 
\ $\overline{\Upsilon}^*$.
\end{proof}
\begin{corollary}\label{T3.11} Assume that \ $\fnouse$ \ is a weak mouse with \ $n(\fnouse)=\overn$. \ Then \ $\fnouse$
\ realizes \ $\Upsilon$ \  if and only if \ $\overline{\core(\fnouse)}$ \ realizes \  
$\overline{\Upsilon}^*$.
\end{corollary}
\begin{proof} Let \ $\alpha$ \ be the ordinal such that the mouse iterate \ $\core(\fnouse)_\alpha$ \ is such that \
$\core(\fnouse)_\alpha=\nouse$. \ Lemma~\ref{lowrho} of this paper and 
Lemma~{\deftwo}  of \cite{Part1} imply that \ $\alpha$ \ is a multiple of \ $\omega^\omega$. \ The result
now follows from Lemma~\ref{T3.10}. 
\end{proof}
\begin{corollary} $\overcore$ \ realizes the type \ $\overline{\Upsilon}^*$.
\end{corollary}
\begin{proof} Since \ $\mouse$ \ is a weak mouse realizing  \ $\Upsilon$, \ Corollary~\ref{T3.11} implies
the desired conclusion. 
\end{proof}
\begin{lemma}\label{applesoranges} Let \ $\nouse$ \ be an iterable pure premouse for which \ $\mouse$ \ is a proper
initial segment.  Let \ $\theta$ \ be any
ordinal which is a multiple of \ $\omega^\omega$ \ such that the premouse iterate \ $\nouse_\theta$ \ and the mouse
iterate \ $\mouse_\theta$ \ are comparable. Then 
\be
\item $\mouse_{\theta}$ \ is a proper initial segment of \ $\nouse_\theta$,
\item $\mouse_{\theta}$ \ realizes \ $\Upsilon$ \ and  \ $(\mouse_{\theta})^\gamma$ \ does not realize 
\ $\Upsilon$ \ for each \ $\gamma < \widehat{\OR}^{\mouse_{\theta}}$,
\item $\mouse_{\theta}$ \ is a weak mouse.
\ee
\end{lemma}
\begin{proof} Let \ $\pi_{0\theta}\colon\nouse\mapsigma{1}\nouse_\theta$ \ be the premouse embedding. Since \
$\mouse\in\nouse$, \ the proof of Theorem~{\intismi} of \cite{Part1} shows that the mouse iterate \
$\mouse_{\theta}$ \ is an initial segment of \ $\pi_{0\theta}(\mouse)$. \ Thus, \ $\mouse_{\theta}$ \ is a proper
initial segment of \ $\nouse_\theta$ \ and therefore, assertion (1) holds. Theorem~\ref{T3.10} implies that \
$\mouse_\theta$ \ realizes the type \ $\Upsilon$ \ and,  because \ $\mouse\in\nouse$, \
we also have that \ $\pi_{0\theta}(\mouse)$ \  realizes the type \ $\Upsilon$. 
\begin{clam} $(\pi_{0\theta}(\mouse))^\gamma$ \ does not realize 
\ $\Upsilon$ \ for any \ $\gamma < \widehat{\OR}^{\pi_{0\theta}(\mouse)}$. 
\end{clam}
\begin{proof}[Proof of Claim] Suppose, for a contradiction, that some \ $\gamma <
\widehat{\OR}^{\pi_{0\theta}(\mouse)}$
\ is such that \ $(\pi_{0\theta}(\mouse))^\gamma$ \ does realize the type \ $\Upsilon$. \ Thus, for some \
$b\in(\pi_{0\theta}(\mouse))^\gamma$ \ we have that 
\[\varphi\in\Upsilon \iff (\pi_{0\theta}(\mouse))^\gamma\models\varphi(b).\]
It follows that \ $\Upsilon$ \ can be interpreted as a real in \ $\nouse_\theta$. \ Hence, \ $\Upsilon$ \ 
is a  `real'  in \ $\nouse$.\footnote{We are not assuming that $\R^\nouse=\R$.}
Therefore, the assertion
\[\left(\exists\gamma \in
\widehat{\OR}^{\pi_{0\theta}(\mouse)}\right)\left(\exists b\in\pi_{0\theta}(\mouse)^\gamma\right)  
\left(\forall\varphi\in\Upsilon\right)
\left[\pi_{0\theta}(\mouse)^\gamma\models\varphi(b)\right].\]
can be expressed as a \ $\Sigma_1$ \ statement, in the parameters \
$\pi_{0\theta}(\mouse)$
\ and the `real'
\
$\Upsilon$, \ which is true in \ $\nouse_\theta$. \ Since \ $\pi_{0\theta}\colon\nouse\mapsigma{1}\nouse_\theta$, \
we conclude that this \ $\Sigma_1$ \ statement, in the parameters \
$\mouse$ \ and the `real' \ $\Upsilon$, \ is true in \ $\nouse$. \ Hence,
\[\left(\exists\gamma \in
\widehat{\OR}^{\mouse}\right)\left(\exists b\in\mouse^\gamma\right)  
\left(\forall\varphi\in\Upsilon\right)
\left[\mouse^\gamma\models\varphi(b)\right].\]
Therefore, \  $\Upsilon$ \ does not witness that \ $\mouse$ \ is weak. This contradiction completes the proof of the
claim.
\end{proof}
Since \ $\mouse_\theta$ \ and \ $\pi_{0\theta}(\mouse)$ \ both realize the type \ $\Upsilon$, \ and because \ $\mouse_\theta$
\ is an initial segment of \ $\pi_{0\theta}(\mouse)$, \ the Claim implies that \
$\mouse_\theta=\pi_{0\theta}(\mouse)$. \ Assertions (2) and (3) follow.
\end{proof}
\begin{lemma}\label{typecomparison1} Let \ $\kouse$ \ be a mouse such that \ $n(\kouse)=\overn$ \ and \
$\R^{\kouse}=\R^{\mouse}$. \  Suppose that \ $\overkouse$ \ realizes \ $\overline{\Upsilon}^*$ \ and that \ $\sigma\colon
\overkouse\mapsigma{1}\overcore$. \ Then \
$\overkouse=\overcore$.
\end{lemma}
\begin{proof} Let \ $\nouse$ \ be an iterable pure premouse for which \ $\mouse$ \ is a proper
initial segment. \ Let \ $\theta$ \ be an ordinal which is a multiple of \ $\omega^\omega$ \ such that the premouse
iterate \ $\nouse_\theta$ \ and the mouse iterates
\ $\kouse_\theta$, \ $\core_\theta$  and \ $\mouse_\theta$ \ are all comparable.\footnote{By combining the proofs of
the Lemma 2.25 of \cite{Crcm} and Lemma 2.23 of \cite{Cfsrm}, one can prove that such a $\theta$
exists.} Thus, \ $\core_\theta=\mouse_\theta$. \ Lemma~\ref{T3.10} implies that \ $\kouse_\theta$ \ realizes \
$\Upsilon$. \  Since \ $\sigma\colon\overkouse\mapsigma{1}\overcore$, \ it follows that \ $\kouse_\theta$ \ must be
an initial segment of \
$\core_\theta=\mouse_\theta$. \ Lemma~\ref{applesoranges}(2) implies that \ $\kouse_\theta=\core_\theta$. \
Therefore, \ $\core(\kouse)=\core$. \ We show that \ $\overkouse=\overcore$. \ For a contradiction, suppose
that \ $\overkouse\ne\overcore$. \ Then there is a premouse iterate, say \ $\overcore_\alpha$ \ with \ $\alpha>0$, \
such that \ $\overcore_\alpha=\overkouse$ \ by Theorem~{\relationtwo} of \cite{Part1}.  Let \
$\pi_{0\alpha}\colon\overcore\mapsigma{1}\overcore_\alpha$
\ be the premouse embedding. Because \ $\sigma\colon\overkouse\mapsigma{1}\overcore$, \ the restriction \
$\sigma\colon \kappa^{\overkouse}\to \kappa^{\overcore}$ \ is an order preserving embedding. Therefore, \
$\kappa^{\overkouse}\le \kappa^{\overcore}$. \ Because \ $\pi_{0\alpha}\colon\overcore\mapsigma{1}\overkouse$ \ is a
premouse embedding,
we obtain that \ $\kappa^{\overkouse}=\pi_{0\alpha}(\kappa^{\overcore})>\kappa^{\overcore}$. \ Contradiction.
\end{proof}
A similar argument proves our next lemma. Recall Definition~{\Extendible} of \cite{Part1}.
\begin{lemma}\label{typecomparison2} Let \ $\kouse$ \ be a mouse such that \ $n(\kouse)=\overn$ \ and \
$\R^{\kouse}=\R^{\mouse}$.
\  Suppose that  \ $\overkouse$ \ realizes \ $\overline{\Upsilon}^*$ \ and that \
$\sigma\colon \F^{\overkouse}\to\F^{\overcore}$ \ is \ $\ul{E}_\overn$--extendible. \ Then \
$\overkouse=\overcore$.
\end{lemma}

We will now focus on the structure of \ $\overcore$. \ We shall show that \ $\overcore$ \ is
the union of a canonical sequence of substructures of \ $\overcore$. \ We shall use the type 
\ $\overline{\Upsilon}^*$ \ to construct this sequence which
will be used in  our proof of Theorem~\ref{newthm}. Recall that \ $m=\overn+k$, \ where \
$m=m(\mouse)$ \ and \ $\overn=n(\mouse)$. \ It follows that \ $\rho_\core^m=1$ \ and so,
\ $\rho_{\overcore}^k=1$. \ Also recall the convention that \ $\overline{C}$ \ denotes the domain of \
$\overcore$. \ The construction of the sequence of substructures is divided into three separate cases;
namely, (1) \ $k>1$, \ (2) \ $k=1$ \ and \ $\widehat{\OR}^{\overcore}$ \ is a limit ordinal, and  (3) \ $k=1$ \ and \
$\widehat{\OR}^{\overcore}$ \ is a successor ordinal. Our next lemma will be used to construct
the desired sequence of substructures in the case when \ $k>1$. \ First recall that Corollary~{\skolem} of
\cite{Part1} implies that there is a \ $\Sigma_{k-1}$ \ Skolem function for \ $\overcore$ \ which is \
$\Sigma_{k-1}(\overcore,\{q,p^{}_{\overcore}\})$ \ for some \ $q\in \overline{C}$.

\begin{lemma}\label{hull} Suppose \ $k\ge 2$ \ and hence, \ $\rho_{\overcore}^{k-1}>1$. \ 
Let \ $q\in \overline{C}$ \ be such that there is a \ $\Sigma_{k-1}$ \ Skolem function for \ $\overcore$ \ which is \
$\Sigma_{k-1}(\overcore,\{q,p^{}_{\overcore}\})$. \ Let \ $\H=\Hull_{k-1}^{\overcore}(\R^{\overcore}\cup
\{q,p^{}_{\overcore}\})$. \ Then
\be
\item $\H\prec_{k-1}\overcore$
\item $\H$ \ does not realize \  $\overline{\Upsilon}^*$
\item $\pow(\R)\cap\boldface{\Sigma}{\omega}(\H)\subseteq\overline{C}$.
\ee
\end{lemma}
\begin{proof} Since \ $\H=\Hull_{k-1}^{\overcore}(\R^{\overcore}\cup
\{q,p^{}_{\overcore}\})$ \ and  there is a \ $\Sigma_{k-1}$ \ Skolem function for \ $\overcore$ \ which is \
$\Sigma_{k-1}(\overcore,\{q,p^{}_{\overcore}\})$, \ we have that \ $\H\prec_{k-1}\overcore$. \ To prove (2),
let \ $\A$ \ be the transitive collapse of \ $\H$ \ and let \ $\pi\colon\A\to\H$ \ be the inverse of the collapse
map. Because \ $\overcore\models\T^\overn$, it follows that \ $\A\models\T^\overn$. \ Lemma~{\EOEL} of
\cite{Part1} implies that 
\be
\item[(a)] there is an acceptable premouse \ $\kouse$ \ such that \ $\kouse^\overn=\A$, 
\item[(b)] there is a map \ $\sigma\colon\kouse\mapsigma{m-1}\core$ \ where \ $m-1\ge\overn+1$.
\ee
Arguing as in the proof of Lemma 2.29 of \cite{Cfsrm}, we have the \ $\kouse$ \ is a mouse and \
$n(\kouse)=\overn$. \ Let \ $\overkouse=\kouse^\overn$ \ and note that \ $\R^{\overcore}=\R^{\overkouse}$. \ Since \
$k>1$ \ we conclude that \ $\pi\colon\overkouse\mapsigma{1}\overcore$.
\begin{clam} There exists an ordinal \ $\theta$ \ such that the mouse iterate \ $\kouse_\theta$ \ is a proper
initial segment of the mouse iterate \ $\core_\theta$.
\end{clam}
\begin{proof}[Proof of Claim] Lemma 2.23 of \cite{Cfsrm} asserts that there is an ordinal \ $\theta$ \ such that
the mouse iterates \ $\kouse_\theta$ \ and \ $\core_\theta$  \ are comparable (see \cite[Definition 2.23]{Crcm}).
Since \ $\pi\colon\overkouse\to\overcore$ \ is \ $\ul{E}_\overn$--extendible, Theorem~{\EOEL} in \cite{Part1} and (the
proof of) Theorem 2.28 in \cite{Cfsrm} imply that we must have either \ $\kouse_\theta=\core_\theta$ \ or \
$\kouse_\theta\in\core_\theta$. \ Suppose, for a contradiction, that \ $\kouse_\theta=\core_\theta$. \ Since \
$\overkouse=\Hull_{k-1}^{\overkouse}(\R^{\overkouse}\cup
\{\pi^{-1}(q),\pi^{-1}(p^{}_{\overcore})\})$ \ and there is a \ $\Sigma_{k-1}$ \ Skolem function for \ $\overkouse$
\ which is \ $\Sigma_{k-1}(\overkouse,\{\pi^{-1}(q),\pi^{-1}(p^{}_{\overcore})\})$, \ it follows that there is a \
$\boldface{\Sigma}{k-1}(\overkouse)$  \ set of reals \ $D$ \ such that \ $D\notin\overkouse$. \ In addition, we
see that \ $D$ \ is also \ $\boldface{\Sigma}{k-1}(\overcore)$. \ Since \ $\rho_{\overcore}^{k-1}>1$, \ we
conclude that \ $D\in\overcore$. \ Hence, by Lemma 2.19 of \cite{Cfsrm} we have that \ $D\in\core_\theta$. \ Because
 \ $\kouse_\theta=\core_\theta$, \ it follows that \ $D\in\kouse_\theta$. \ Again by by Lemma 2.19 of \cite{Cfsrm},
we conclude that \ $D\in\overkouse$. \ This contradiction completes the proof of the claim.
\end{proof}
Since \ $\pi\colon\overkouse\mapsigma{1}\overcore$, \ Lemma~\ref{typecomparison1} and the above Claim now imply that \
$\overkouse$ \ does not realize the type \ $\overline{\Upsilon}^*$. \ Therefore,
\ $\H$ \ does not realize the type \ $\overline{\Upsilon}^*$ \ and this completes the proof of (2). To
establish (3), we observe that the above Claim and Lemma 2.19 of \cite{Cfsrm} also imply that \
$\pow(\R)\cap\boldface{\Sigma}{\omega}(\H)\subseteq\overline{C}$.
\end{proof}

\begin{remark}\label{allreals} For the remainder of this section we shall assume that  \ $\mouse$ \ is a real mouse. \
Thus, \ $\R^{\mouse}=\R^{\overcore}=\R$.
\end{remark}
The concept of a good covering will be used to construct scales in our proof of Theorem
\ref{goodcovering} below. This concept will also be used to evaluate the complexity of these scales. We
shall now state the definition of a good covering given in \cite[Definition~{\goodcovdef}]{Part1}.
Given a structure \ $\H$ \ we shall let \ $H$ \ denote the domain of this structure.
\begin{definition}\label{goodcoveringdef} Let \ $\nouse =
(N,\in,\R,c_1,c_2,\dots,c_m,A_1,A_2,\dots,A_n)$ \ 
be a transitive model of \ $\Rplus$. \ 
Suppose that \ $\langle\mathcal{H}^i : i\in\omega\rangle$ \ is a sequence of substructures of \
$\nouse$ \ such that \ $\R\subseteq H^0\subseteq H^1\subseteq\cdots \subseteq\bigcup\limits_{i=0}^\omega
H^i=N$.
\ We shall say that \ $\langle\mathcal{H}^i : i\in\omega\rangle$ \ is a {\it good covering\/} of \ $\nouse$ \ if
for each \
$i\in\omega$,
\be
\item $\mathcal{H}^i\prec_0\mathcal{H}^{i+1}$
\item $\mathcal{H}^i$ \ is in \ $\boldface{\Sigma}{\omega}(\mathcal{H}^{i+1})$
\item $\pow(\R)\cap\boldface{\Sigma}{\omega}(\mathcal{H}^i)\subseteq N$
\item there exists a \ $\boldface{\Sigma}{\omega}(\mathcal{H}^{e(i)})$ \ function \
$f\colon\R\maps{onto}H^i$,
\ee
for some fixed \ $e\in\R$ \ where \ $e(j)\ge j$ \ for all \ $j\in\omega$.
\end{definition}
 
The proof of Theorem~\ref{goodcovering} shows that \ $\overcore$ \ has a good covering in all but one case. However, 
even in this case we can construct a covering of \ $\overcore$ \ that satisfies 
conditions (1)-(3) of the above definition. This fact motivates the following definition.

\begin{definition}\label{coveringdef} Let \ $\nouse =
(N,\in,\R,c_1,c_2,\dots,c_m,A_1,A_2,\dots,A_n)$ \ 
be a transitive model of \ $\Rplus$. \ 
Suppose that \ $\langle\mathcal{H}^i : i\in\omega\rangle$ \ is a sequence of substructures of \
$\nouse$ \ such that \ $\R\subseteq H^0\subseteq H^1\subseteq\cdots \subseteq\bigcup\limits_{i=0}^\omega
H^i=N$. \ We shall say that \ $\langle\mathcal{H}^i : i\in\omega\rangle$ \ is a {\it suitable covering\/} of \ $\nouse$ \ if
for each \ $i\in\omega$ \ conditions (1), (2) and (3) in Definition~\ref{goodcoveringdef} hold.
\end{definition}

If the structure \ $\overcore$ \ has a good covering and satisfies \ $\AD$, \ then \ $\overcore$ \ posses
some interesting definability properties which will be used to show that a certain Moschovakis scale exists and is definable
over \ $\overcore$.

\begin{definition} Let \ $\Lng = \Lng_n\cup\{B_1,\dots,B_i,\dots\}$ \ 
where the \ $B_i$'s \ are new predicate symbols \ for \ $1\le i<\omega$. \ A quantifier is {\it
bounded\/} in \ $\Lng$ \  if it has the form \ $\exists u\in v$, \ $\forall u\in v$, \ $\exists
u\in B_i$, \ or \ $\forall u\in B_i$. \  A formula \ $\varphi$ \ in \ $\Lng$
\ is said to be in \
$\Sigma_0^0$ \ if all the quantifiers in
\ $\varphi$ \ are  bounded. 
\end{definition}

\begin{definition}\label{stuff}
Let \ $\nouse$ \ be a transitive model of \ $\Rplus$ \ and let  
\ $\langle\H_i : i\in\omega\rangle$ \ be a sequence of substructures of \ $\nouse$. \ For \ $k\in\omega$, \ we say that \
$P\subseteq N\times\omega^k$ \ is in \
$\boldface{\Sigma}{0}^\omega(\nouse,H_0,H_1,\dots)$
\ if there is a \ $\Sigma_0^0$ \ formula \
$\varphi(u,v_1^\prime,\dots,v_k^\prime,v_1,\dots,v_m,B_1^\prime,\dots,B_k^\prime,B_1,\dots,B_m)$ \ (allowing
constants from
\ $N$) \ such that \ $(j_1,\dots,j_k,a)\in P$ \ if and only if
\[\forall i_1\exists i_2\cdots Q i_m\,(\nouse,H_0,H_1,\dots)\models
\varphi(a,j_1,\dots,j_k,i_1,\dots,i_m,H_{j_1},\dots,H_{j_k},H_{i_1},\dots,H_{i_m})\]
where \ $Q i_m$ \ is either \ $\forall i_m$ \ or \ $\exists i_m$ \ depending on the parity of \ $m$.
\end{definition}

Let \ $\nouse$ \ and \ $\langle\H_i : i\in\omega\rangle$ \ be as in Definition~\ref{stuff}. Let \ $\varphi$ \ be 
a \ $\Sigma_0^0$ \ formula, say \ $\varphi=\varphi(u,v_1,\dots,v_{m},B_1,\dots,B_m)$. \ For \ $s\in\omega$ \
and \ $\vec\iota=\langle i_1,\dots,i_m\rangle\in{{^{m}}\omega}$ \ we shall write
\[\varphi_{\vec\iota}^s=\{a\in H_s :
(\nouse,H_0,H_1,\dots,H_{i_m})\models \varphi(a,i_1,\dots,i_{m},H_{i_1},\dots,H_{i_m})\}.\] In addition,
given a prewellordering \ $\le$ \ on an arbitrary set, we say that \ $\thickapprox$ \ is the {\it equivalence relation
derived\/} from \ $\le$ \ when \ $x\thickapprox y \Leftrightarrow x\le y\land y\le x$.

\begin{definition} 
Let \ $\nouse$ \ be a transitive model of \ $\Rplus$, \ $\langle\H_i : i\in\omega\rangle$ \ 
a sequence of substructures of \ $\nouse$, \ and \
$\varphi=\varphi(u,v_1,\dots,v_{m},B_1,\dots,B_m)$ \ in \ $\Sigma_0^0$. \ Given \ $s\in\omega$ \
and \ $Q\in\boldface{\Sigma}{\omega}(\H_s)$, \ suppose that \ $Q$ \ has a prewellordering \
$\le\ \in\boldface{\Sigma}{\omega}(\H_s)$ \ and let \ $\thickapprox$ \ be the equivalence relation derived
from  \ $\le$. \ 
If for all \ $\vec\iota=\langle i_1,\dots,i_m\rangle\in{{^{m}}\omega}$ \ 
\be
\item $\varphi_{\vec\iota}^s\subseteq Q$
\item $x\thickapprox y \implies [\varphi_{\vec\iota}^s(x) \Leftrightarrow \varphi_{\vec\iota}^s(y)]$, \
for all \
$x, y \in Q$
\ee
then \ $\varphi$ \ is said to be \ $\thickapprox^{\H_s}$--{\it invariant}.
\end{definition} 
In the above definition, note that if \ $\le$ \ is a well-ordering on \ $Q$, \ then 
condition (2) holds trivially. The following is Theorem 3.31 of \cite{Cfsrm}.

\begin{theorem}\label{ADdefinable} Let \ $\nouse$ \ be a transitive model of \ $\Rplus + \AD$, \ containing all the reals, \
and let \ $\langle\H_i : i\in\omega\rangle$ \ be a good covering of \ $\nouse$. \ Given \ $s\in\omega$, \ 
suppose that \ $P\subseteq H_s$ \ is in \ $\boldface{\Sigma}{0}^\omega(\nouse,H_0,H_1,\dots)$ \ as witnessed by
\ $\varphi\in\Sigma_0^0$. \ If \ $\varphi$ \ is \ $\thickapprox^{\H_s}$--{\it invariant}, \  then \ $P$ \ is in \
$\boldface{\Sigma}{\omega}(\H_{e(s)})$.
\end{theorem}

\begin{theorem}[$\ZF+\DC$]\label{goodcovering} The structure \ $\overcore$ \ has a suitable covering \ $\langle\mathcal{H}^i
: i\in\omega\rangle$ \ with the following property: If \ $P\subseteq\R$ \ is \ $\boldface{\Sigma}{k}(\H^{i_0})$ \ for some
fixed \ $i_0\in\omega$, \ then \ $P$ \ has a closed game representation and, assuming \ $\overcore\models\AD$, \ the set \ 
$P$ \ has a scale that is \ $\boldface{\Sigma}{k}(\overcore)$. \ Moreover, there is a partial \
$\boldface{\Sigma}{k}(\overcore)$ \ map of \ $\R$ \ onto \ $\overC$, \ the domain of \ $\overcore$.
\end{theorem}

\begin{rmk}One who is familiar with Steel's analysis of scales in \ $\Lr$ \ and his proof of Theorem 3.7 in
\cite{steel} will see that our proof of Theorem~\ref{goodcovering} is motivated and inspired by Steel's solution to the problem of finding scales at the end of a gap in \ $\Lr$. 
\end{rmk}

The proof of  Theorem~\ref{goodcovering} contains a series of lemmas, beginning with Lemma~\ref{lemma_k>1} and ending with Lemma~\ref{finally}.

\begin{proof}[Proof of Theorem~\ref{goodcovering}]
We will first standardize the parameter in \ $\overC$ \ that 
realizes the type \ $\overline{\Upsilon}^*$. \ Let \ $b'\in \overC$ \ be such that \
$\overcore\models\theta(b')$ \ for all \ $\theta\in\overline{\Upsilon}^*$. \ Corollary 1.8 of \cite{Crcm} implies
that there is a uniformly \ $\Sigma_1(\overcore)$ \ function
\ $f^{\overcore}\colon[\OR^{\overcore}]^{<\omega}\times\R\maps{onto}\overC$. Let \
$w_0\in\R$, \ and \ $\G'$ \ be a finite subset of \ $\OR^{\overcore}$ \ such that \
$b'=f^{\overcore}(\G',w_0)$. \ Set
\[\Gamma=\Upsilon^\overn_k(\seq{\G',w_0 },\overcore)\]
(see Definition~\ref{type}). Clearly, \ $\overline{\Upsilon}^*\subseteq\Gamma$.
\ Let \ $\G\le_{\BK}\G'$ \ be the $\le_{\BK}$--least \ $\G$ \ such that
\[\Gamma=\Upsilon^\overn_k(\seq{\G,w_0 },\overcore).\]
Let \ $<_\Gamma$ \ be a fixed order of \ $\Gamma$ \ in order type \ $\omega$.

Using \ $\G$ \ and \ $w_0$ \ we can now define our desired sequence of substructures. The proof is broken into three cases: 
\be
\item[(1)] $k>1$
\item[(2)] $k=1$ and ${\widehat{\OR}^{\overcore}}$ is a limit ordinal
\item[(3)] $k=1$ and ${\widehat{\OR}^{\overcore}}$ is a successor ordinal.
\ee

{\bf Case 1: $\mathbf{k>1}$.}  \ Let \ $S'\colon\R\times\overC\to\overC$ \ be a \ $\Sigma_{k-1}$ \
Skolem function for \
$\overcore$. \ Corollary~{\skolem} of \cite{Part1} implies that \ $S'$ \ exists. Let \ $\theta_{\textup{sk}}\in
\Sigma_{k-1}$, \ $w_1\in\R$, \ and \ $F'$ \ be a finite subset of \
$\OR^{\overcore}$ \ such that
\[S'(x,u)=v\iff \overcore\models \theta_{\textup{sk}}(x,u,v,w_1,F').\]
We shall standardize the parameter \ $F'$ \ that will give such a Skolem function.
Let \ $F\le_{\BK}F'$ \ be the $\le_{\BK}$--least \ $F$ \ such that \
$S\colon\R\times\overC\to\overC$ \ defined by
\[S(x,u)=v\iff \overcore\models \theta_{\textup{sk}}(x,u,v,w_1,F)\]
is a \ $\Sigma_{k-1}$ \ Skolem function for \ $\overcore$. 

Using \ $F$, \ $\G$, \ $w_0$ \ and Lemma~\ref{hull}, we will define a canonical sequence \
$\seq{\H^i : i\in\omega}$ \ of \ $\Sigma_{k-1}$ \ hulls. At the same time we shall
also define a sequence \ $\seq{\theta_i : i\in\omega}$ \ of \ $\Sigma_k$ \ formula in
\ $\Gamma$, \ and a sequence \ $\seq{Q_i : i\in\omega}$ \ where each \ $Q_{i+1}$ \ is a
finite subset of \ $\OR^{\overcore}$.

First, define \ $Q_0=\seq{\rho^{}_{\,\overcore}, p^{}_{\,\overcore}}$. \ Assuming that
\
$Q_0, Q_1, \dots, Q_i$ \  have been defined, we shall define \ $\H^i$, \ $\theta_i$ \ and \
$Q_{i+1}$ \ as follows:

\smallskip
\noindent{\sl Definition of \ $\H^i$\/}: \ Define
\ $\H^i=\Hull_{k-1}^{\overcore}(\R\cup\{F, \G, Q_0, Q_1, \dots, Q_i\})$. 
 
By our choice of \ $F$, \ $\G$, \ $Q_0$ \ and \ $\Gamma$, \ Lemma~\ref{hull} implies
\be
\item $\H^i\prec_{k-1}\overcore$
\item $\H^i$ \ does not realize \  $\Gamma$
\item $\pow(\R)\cap\boldface{\Sigma}{\omega}(\H^i)\subseteq\overline{C}$.
\ee

\smallskip
\noindent{\sl Definition of \ $\theta_i$\/}: \ Since \ $\H^i$ \ does not realize \ $\Gamma$, \ define \
$\theta_i$ \ to be the $<_\Gamma$\,--least 
\ $\Sigma_k$ \ formula \ $\theta(v)\in\Gamma$ \ such that 
\ $\H^i\not\models\theta(\seq{\G,w_0 })$.

\smallskip
\noindent{\sl Definition of \ $Q_{i+1}$\/}: \ Because \ $\theta_i$ \ is \ $\Sigma_k$ \ and
because every element in \ $\overC$
\ is \ $\Sigma_1$ \ definable from a real in \ $\R$ \ together with a finite subset of \
$\OR^{\overcore}$ (see
\cite[Corollary 1.8]{Crcm}), there is a finite \ $Q\subseteq\OR^{\overcore}$ \ such that 
\begin{equation}
\Hull_{k-1}^{\overcore}(\R\cup\{F, \G, Q_0, Q_1, \dots, Q_i,
Q\})\models\theta_i(\seq{\G,w_0}).
\label{maketrue}\end{equation} 
Let \ $Q_{i+1}$ \ be the $\le_{\BK}$--least such \ $Q$ \ satisfying (\ref{maketrue}).

\smallskip
This completes the definition of the canonical sequences \ $\seq{\H^i :
i\in\omega}$, \ $\seq{\theta_i :
i\in\omega}$ \ and \ $\seq{Q_i :
i\in\omega}$. \ 
We now record some properties of these sequences.
\begin{lemma}\label{lemma_k>1} Let \ $\seq{\H^i :
i\in\omega}$, \ $\seq{\theta_i :
i\in\omega}$ \ and \ $\seq{Q_i :
i\in\omega}$ \ be as above. Then \ $\bigcup\limits_{i\in\omega}\H^i=\overcore$ \ and for all \ $i\in\omega$ \ 
\be
\item[(a)] $\H^i\prec_{k-1}\H^{i+1}\prec_{k-1}\overcore$
\item[(b)] $\mathcal{H}^i$ \ is in \ $\boldface{\Sigma}{k-1}(\mathcal{H}^{i+1})$
\item[(c)] $\pow(\R)\cap\boldface{\Sigma}{\omega}(\H^i)\subseteq\overline{C}$
\item[(d)] there exists a \ $\boldface{\Sigma}{k-1}(\mathcal{H}^{i})$ \ partial function \ $f\colon\R\maps{onto}H^i$
\item[(e)] $\H^i\not\models\theta_i(\seq{\G,w_0})$  \ and \  $\H^{i+1}\models\theta_i(\seq{\G,w_0})$.
\ee
\end{lemma}
\begin{proof} For each \ $i\in\omega$, \ items (a)--(e) are clear. To prove that \
$\bigcup\limits_{i\in\omega}\H^i=\overcore$,
\ let
\ $\H=\bigcup\limits_{i\in\omega}\H^i$. \ Thus,
\ $\H\prec_{k-1}\overcore$. \ Arguing as in the proof of Lemma 2.29 of \cite{Cfsrm}, there is a mouse \ $\kouse$ \ such
that \ $n(\kouse)=n(\core)=\overn$ \ and there is an  isomorphism \ $\pi\colon\H\to\overkouse$ \ where \ $\pi$ \
is the transitive collapse of \ $\H$. \ Therefore, \ $\overkouse$ \ realizes the type \ $\overline{\Upsilon}^*$. \
Since \
$\pi^{-1}\colon\overkouse\mapsigma{1}\overcore$, \ Lemma~\ref{typecomparison1} implies that \ $\overkouse=\overcore$.
\ Thus, \ $\pi\colon\H\to\overcore$ \ is an (onto) isomorphism. The following two claims establish that \
$\H=\overcore$.
\begin{clm} Assume that \ $\pi(F)=F$, \ $\pi(\G)=\G$ \ and \ $\pi(Q_i)=Q_i$ \ for all \ $i\in\omega$. \ Then \
$\H=\overcore$. 
\end{clm}
\begin{proof}[Proof of Claim 1] Clearly, \ $\pi(x)=x$ \ for all \ $x\in\R$. \ Since \
\[\H=\bigcup\limits_{i\in\omega}\H^i=\Hull_{k-1}^{\H}(\R\cup
\{F,\G,Q_0,\dots,Q_i,\dots\}),\] it follows that 
\[\overcore=\Hull_{k-1}^{\overcore}(\R\cup
\{\pi(F),\pi(\G),\pi(Q_0),\dots,\pi(Q_i),\dots\}).\]
Therefore, our assumption implies that \ $\H=\overcore$.
\end{proof}

\begin{clm} $\pi(F)=F$, \ $\pi(\G)=\G$ \ and \ $\pi(Q_i)=Q_i$ \ for all \ $i\in\omega$. 
\end{clm}
\begin{proof}[Proof of Claim 2] To see that \ $\pi(F)=F$, \ note that \ $\theta_{\textup{sk}}$, \ in parameters
\ $w_1\in\R$ \ and \ $F$, \ defines a \ $\Sigma_{k-1}$ \ Skolem function over \ $\H$. \ Therefore, \ 
$\theta_{\textup{sk}}$, \ in parameters
\ $w_1$ \ and \ $\pi(F)$, \ defines a \ $\Sigma_{k-1}$ \ Skolem function over \ $\overcore$. \ But \ $\pi(F)\le_{\BK}F$
\ and \ $F$ \ is the $\le_{\BK}$--least such set of ordinals. Hence, \ $\pi(F)=F$. \ Similarly, \ one can see that \
$\seq{\pi(\G),w_0}$ \ realizes the type \ $\Gamma$ \ in \ $\overcore$. \ Since \ $\pi(\G)\le_{\BK}\G$, \ the
$\le_{\BK}$--minimality of \ $\G$ \ implies that \ $\pi(\G)=\G$.

Finally, we show that \ $\pi(Q_i)=Q_i$ \ by induction on \ $i$. \ We first show that \
$\pi(Q_0)=Q_0$, \ where \ $Q_0=\seq{\rho^{}_{\,\overcore}, p^{}_{\,\overcore}}$. \ Lemma~{\soundnplusone} of \cite{Part1}
implies that \
$\overcore$ \ is sound \ (see \cite[Definition~{\sounddef}]{Part1}). Therefore, we have that \
$h_{\overcore}^{p^{}_{\,\overcore}}\colon\R\times\omega\rho^{}_{\,\overcore} \maps{onto} \overC$ \ (see \cite[Definition
{\skolemdef}]{Part1}). This assertion can be expressed as a \ $\Pi_2$ \ statement, in the parameters  \
$p^{}_{\,\overcore}$ \ and \ $\rho^{}_{\,\overcore}$, \ which is true in \ $\overcore$. \ Since \
$\H\prec_{1}\overcore$, \ it follows that \ $h_{\H}^{p^{}_{\,\overcore}}\colon
\R\times\omega\rho^{}_{\,\overcore} \maps{onto} H$, \ where \ $H$ \ is the domain of the structure \
$\H$. \ Thus, there is a \ $D\subseteq \R\times \omega\rho^{}_{\,\overcore}$ \ which is \
$\Sigma_1(\H,\{p^{}_{\,\overcore}\})$ \ and \ $D\notin\H$. \ Since \ $\pi\colon\H\to\overcore$ \ is an isomorphism, we
conclude there must be a \ $D'\subseteq \R\times\omega\pi(\rho^{}_{\,\overcore})$ \ which is \
$\Sigma_1(\overcore,\{\pi(p^{}_{\,\overcore})\})$ \ and \ $D'\notin\overcore$. \ But \
$\pi(\rho^{}_{\,\overcore})\le\rho^{}_{\,\overcore}$ \ and \ $\rho^{}_{\,\overcore}$ \ is the least ordinal such that
there is a \ $\boldface{\Sigma}{1}(\overcore)$ \ subset of \ $\R\times \omega\rho^{}_{\,\overcore}$ \ which
is not in \ $\overcore$. \ Therefore, \ $\pi(\rho^{}_{\,\overcore})=\rho^{}_{\,\overcore}$. \ Similarly, one can show
that \ $\pi(p^{}_{\,\overcore})=p^{}_{\,\overcore}$. \ Hence, \ $\pi(Q_0)=Q_0$. \ Assume that ($*$) 
$\pi(Q_j)=Q_j$ \ for all \ $j\le i$. \ Let \ $\pi\l\H^{i+1}\r$ \ be the isomorphic image of the structure \ $\H^{i+1}$.
\ Clearly, \ $\pi\l\H^{i+1}\r$ \ is a substructure of \ $\overcore$ \ and since
\[\H^{i+1}=\Hull_{k-1}^{\overcore}(\R\cup\{F, \G, Q_0, Q_1, \dots, Q_i, Q_{i+1}\})\models
\theta_i(\seq{\G,w_0 }),\]
$\pi(F)=F$ \ and \ $\pi(\G)=\G$, \ it follows from ($*$) that 
\[\pi\l\H^{i+1}\r=\Hull_{k-1}^{\overcore}(\R\cup\{F, \G, Q_0, Q_1, \dots, Q_i, \pi(Q_{i+1})\})\models
\theta_i(\seq{\G,w_0 }).\]
Since \ $\pi(Q_{i+1})\le_{\BK}Q_{i+1}$ \ and \ $Q_{i+1}$ \ is the $\le_{\BK}$--least such set of ordinals, we must
have that \ $\pi(Q_{i+1})=Q_{i+1}$. \  This completes the proof of Claim 2.
\end{proof}
The proof of the lemma is complete. 
\end{proof}
Letting \ $e(i)=i$ \ for all \ $i\in\omega$, \ it follows that  \ $\langle\H^i : i\in\omega\rangle$ \
is a good covering of \ $\overcore$.
\begin{lemma}\label{onto1} There is a partial \ $\boldface{\Sigma}{k}(\overcore)$ \ map of \ $\R$ \ onto \ $\overC$, \ the
domain of \ $\overcore$.
\end{lemma}
\begin{proof} Because \ $k\ge 2$, \ it follows that \ $\rho_{\overcore}^{k-1}>1$. \ Lemma~{\soundnplusone} of
\cite{Part1} implies that
\ $\overcore$ \ is \ $k$--sound. Since \ $\rho_{\overcore}^{k}=1$, \ Corollary 1.32 of \cite{Cfsrm} states that there is
a partial \ $\boldface{\Sigma}{k}(\overcore)$ \ map of \ $\R$ \ onto \ $\overC$.
\end{proof}
\begin{lemma}\label{closedk>1} Consider the parameters \ $\y,F, \G, Q_0, Q_1, \dots, Q_\ell$ \ (as defined above) where \
$\y$ \ is a fixed real.  Let \ $\theta$ \ be a \ $\Sigma_k$ \ formula in the language \ $\Lng_\overn$ \ augmented by
these parameters. Fix \ $i_0\ge\ell$ \ and let \ $P\subseteq\R$ \ be defined by \ $P(x)$ \ if and only if \
$\H^{i_0}\models \theta(x)$. \ Then \ $P$ \ has a closed game representation and if \ $\overcore\models\AD$, \ then there
is a \ $\boldface{\Sigma}{k}(\overcore)$ \ scale on \ $P$.
\end{lemma}
\begin{proof} We shall construct a closed game representation of \ $P$ \ simple enough to yield a
\ $\boldface{\Sigma}{k}(\overcore)$ \ scale on \ $P$. \ Our closed game representation \ $x\mapsto G_x$ \
of \ $P$ \ will have the following property: If \ $P_n$ \ is defined by
\[P_n(x,u)\iff \textup{$u$ is a winning position for player $\mathbf{I}$ in $G_x$ of length $n$},\]
then \ $P_n$ \ is \ $\boldface{\Sigma}{\omega}(\H^{\max(i_0,n)})$. \ Assuming \ $\overcore\models\AD$, \ Lemma
\ref{poffinmouse} and Lemma~\ref{lemma_k>1}(c) will then imply that the Moschovakis scale \ $\seq{\varphi_i : i \in \omega}$
\ on \ $P$ \ exists. Furthermore, the prewellordering \ $\le_i$ \ of \ $\R$ \ induced by \ $\varphi_i$ \ is also \
$\boldface{\Sigma}{\omega}(\H^{\max(i_0,i)})$. \ Lemma~\ref{lemma_k>1}(c) implies that \ $\le_i\,\in\overC$, the domain of \
$\overcore$, \ for all \ $i\in\omega$. \ Lemma~\ref{onto1} states that there is a partial \ $\boldface{\Sigma}{k}(\overcore)$
\ map of \ $\R$ \ onto \ $\overC$. \ It follows that any countable subset of \ $\overC$ \ is \
$\boldface{\Sigma}{k}(\overcore)$. \ Therefore, \ $\seq{\varphi_i : i \in \omega}$ \ is a \ $\boldface{\Sigma}{k}(\overcore)$
\ scale on \ $P$. 

So, to prove the lemma, it suffices to construct a closed game representation \ $x\mapsto G_x$ \ so that \ $P_n$ \ is
first order definable over the structure \ $\H^{\max(i_0,n)}$ \ for all \ $n\in\omega$. \ In our construction of \ $G_x$ \ 
we shall force player $\mathbf{I}$ to describe the truth in \ $\overcore$ \ about \ $F$, $G$ and the $Q_i$'s. \ To ensure
that each \ $P_n$ \ is \ $\boldface{\Sigma}{\omega}(\H^{\max(i_0,n)})$ \ we impose some restrictions on player $\mathbf{I}$'s
moves. For example, player $\mathbf{I}$ must describe only \ $\Sigma_{k-1}$ \ truths and, to prove that his model is mouse
iterable, player $\mathbf{I}$ must play  `functions' from \ $\F^j$ \ where \ $\F^j=\F^{\H^j}$, \ at each move \ $j$.

Player $\mathbf{I}$'s description of his model is in the language
\[\Lng=\Lng_\overn\cup\{\uF,\uG,\urho,\up\}\cup\{\ux_i,\uQ_i: i\in\omega\}\]
where \ $\overn=n(\core)$. \ If \ $\varphi$ \ is an \ $\Lng$--formula containing no constants \ $\uQ_i$ \ for \ $i>n$,
\ then we shall say that \ $\varphi$ \ has {\it support\/} $n$. 

Let \ $\B_{k-1}$ \ be the set of boolean combinations of \ $\Sigma_{k-1}$ \ formulae in the language \ $\Lng$.
\ Player $\mathbf{I}$ will describe the truth of formulae in \ $\B_{k-1}$. \  We shall use the ``unique $v$'' operator
\ $\iota v$ \ applied to \ $\Sigma_{k-1}$ \ formulae as abbreviations of formulae in \ $\B_{k-1}$. \ For example, let \
$\sigma$, \ $\tau$, \ $\delta$ \ and \ $\varphi$ \ be \ $\Sigma_{k-1}$ \ formulae and let \ $\psi$ \ be a
$\Pi_{k-1}$ \ formula. Then one can check that \ $\varphi(\iota v\sigma(v))$, \ $\psi(\iota v\sigma(v))$ \ and \
$\psi(\iota v\delta(v,\iota u\tau(v)))$ \ can easily be interpreted as formulae in \ $\B_{k-1}$. \ Let \ $T$ \ be 
the \ $\Lng$ \ theory consisting of the sentences
\[\uF\in [\OR]^{<\omega}, \ \uG\in [\OR]^{<\omega}, \ \urho\in\OR, \ \up\in[\OR]^{<\omega}, \ \ux_i\in\ul{\R}, \
\uQ_i\in [\OR]^{<\omega}\]  for each \ $i\in\omega$. \ Clearly, \ $T\subseteq\B_{k-1}$.

In \cite[section~{\realmice}]{Part1} we demonstrated that \
$\overmouse$ \ satisfies the \ $\Lng_\overn$ \ theory \ $\T^\overn$ \ which is defined in \cite[pp.
951-954]{Cfsrm}. We also reviewed in \cite[section~{\coremice}]{Part1} the definition of \ $\core$ \ and the fact that \
$\overcore$ \ satisfies the  theory \ $\T^\overn$. \ Also, as noted in \cite{Cfsrm}, the \ theory \ $\T^\overn$ \ is
axiomatized by a \ $\Pi_2$ \ sentence \ $\forall v\psi^\overn(v)$ \ in the language \ $\Lng_\overn$, \ where \
$\psi^\overn(v)$ \ is \ $\Sigma_1$. \ We shall require player $\mathbf{I}$ to describe a model of the theory \
$\T^\overn$.

For expository reasons, player $\mathbf{I}$ is allowed to play finitely many sentences and finitely many reals in a
single move of  \ $G_x$. \ A typical run of the game has the form
\begin{equation}
\begin{aligned}[c]
{}&{\mathbf{I}}\phantom{{\mathbf{I}}} \qquad T_0,s_0,\f_0,m_0 \qquad \phantom{y_0} \qquad T_1,s_1,\f_1,m_1 \qquad
\phantom{y_1}\quad
\\ {}&{\mathbf{I}}\mathbf{I} \qquad \phantom{T_0,s_0,\f_0,m_0}\qquad y_0 \qquad \phantom{T_1,s_1,\f_1,m_1} \qquad y_1 \quad 
\end{aligned}
\begin{gathered}[c]
{\cdots}
\end{gathered}
\tag*{$G_{x}$}\end{equation}
 where 
\bi
\item $T_j$ \ is a finite set of sentences each of which is in \ $\B_{k-1}$ \ and has support \ $j$ 
\item $s_j\in\R^{<\omega}$, \ $\f_j\in \F^{j}$, \ $j<m_j\in\omega$, \ and \ $y_j\in\R$.
\ei
Given a run  of the game as above, let
\[\seq{x_i:i\in\omega}=s_0{^\frown}\seq{y_0}{^\frown}s_1{^\frown}\seq{y_1}{^\frown}\cdots\]
and let \ $T^*=\bigcup\limits_{j\in\omega}T_j$. \ Let \ $n\colon\B_{k-1}\maps{1-1}\omega$ \ be such that for any \
$\psi\in\B_{k-1}$ \ has support \ $n(\psi)$ \ and has no constants \ $\ux_i$ \ for \ $i\ge n(\psi)$. 

We shall say that the above run is a {\it winning run\/} for player $\mathbf{I}$ if it meets the following 
{\sl closed\/} requirements R-1 to R-9. 
\be
\item[R-1:] $s_0(0)=x$, \ $s_0(1)=w_1$, \ $s_0(2)=w_0$, \ $s_0(3)=\y$.
\item[R-2:]
\be
\item The theory \ $T\cup T^*$ \ is consistent.
\item If \ $\psi$ \ is a sentence in \ $\B_{k-1}$, \ then either \ $\psi\in T_{n(\psi)}$ \ or \ $(\lnot\psi)\in
T_{n(\psi)}$.
\item If \ $\psi\in T_j$, \ then \ $\psi$ \ has support \ $j$ \ and does not contain any constants of the form \
$\ux_i$ \ for \ $i\ge \dom(s_0{^\frown}\seq{y_0}{^\frown}s_1{^\frown}\seq{y_1}{^\frown}\cdots
{^\frown}s_{j-1}{^\frown}\seq{y_{j-1}}{^\frown}s_j)$.
\item If \ $\tau(v)$ \ is \ $\Sigma_{k-1}$ \ and \ $\exists u(u=\iota v\tau(v))$ \ is in \ $T_j$, \ then
\ $\psi^\overn(\iota v\tau(v))$ \ is in \ $T_{j+1}$.
\item $(\ux_j(\un)=\um)$ \ is in \ $T^*$ \ if and only if \ $x_j(n)=m$, \ for all \ $n,m\in\omega$. 
\ee
\ee
\begin{rmk} Requirements R-2(a) and R-2(b) ensure that theory \ $T\cup T^*$ \ has a model and that \ $T^*$ \ is
complete with respect to \ $\B_{k-1}$ \ sentences. Requirement R-2(c) forces player $\mathbf{I}$ to make assertions
only about the reals that have previously been played. Requirement R-2(d) ensures that player $\mathbf{I}$ plays a
model of the theory
\ $\T^\overn$. \ Requirement R-2(e) forces player $\mathbf{I}$ to accurately interpret the reals played. Here \ $\un$
\ and \ $\um$ \ represent canonical representations, in the language \ $\Lng$, \ of the natural numbers \ $n$
\ and \ $m$. 
\end{rmk}

Let \ $\gamma(u,v)$ \ be the \ $\Sigma_{k-1}$ \ formula in the language \ $\Lng_\overn$, \ containing only
the parameters \ $\y, F, \G, Q_0, Q_1, \dots, Q_\ell$, \ so that \ $\theta(u)=\exists v \gamma(u,v)$.
\be
\item[R-3:] For some \ $\Sigma_{k-1}$ \ formula \ $\tau(v)$ \ with support \ $i_0+1$, \ the \ $\B_{k-1}$ \ formula \ 
$\gamma(\ux_0,\iota v\tau(v))$ \ is in \ $T_{i_0+1}$.
\item[R-4:] If \ $\tau(v)$ \ is \ $\Sigma_{k-1}$ \ and the formula \ $(\iota v\tau(v)\in\uR)$ \ is in \ $T_j$, \ then
the  formula \ $(\iota v\tau(v)=\ux_i)$ \ is in \ $T_{j+1}$ \ for some \ $i\in\omega$.
\item[R-5:] If \ $\sigma$ \ and \ $\tau$ \ are \ $\B_{k-1}$ \ formulae of the form \ $(\iota v\phi(v)\in\F_n)$ \ and
\ $(\iota v\psi(v)\in\F_m)$, \ respectively, and if \ $\sigma,\tau$ \ are in \ $T^*$, \ then 
\[\text{$\left(\iota v\phi(v)\,\Eq^{s,t}\,\iota v\psi(v)\right)$  is in \ $T^*$ \ iff \ $\overcore\models
\f_{n(\sigma)}\,\Eq^{s,t}\,\f_{n(\tau)}$}\]
for all \ $s\in(n+m)^n\!\uparrow$ \ and all \ $t\in(n+m)^m\!\uparrow$.
\ee
\begin{rmk} Requirement R-3 ensures that \ $\theta(x)$ \ holds in the model that player $\mathbf{I}$ is constructing.
Requirement R-4 forces player $\mathbf{I}$ to verify that any real he describes is one of the \ $\ux_i$'s. \
Requirement R-5 compels player $\mathbf{I}$ to establish that his model is $\overn$--iterable (see section
{\minmouseiterable} of \cite{Part1} and, in particular, see Theorem~{\criterion} of \cite{Part1}).
\end{rmk}
Our final four requirements ensure that player $\mathbf{I}$ correctly interprets each of the constant symbols
\ $\uF$, $\uG$, $\urho$, $\up$, $\uQ_0, \uQ_1, \dots$ \ in the model he is constructing. To do this, player
$\mathbf{I}$ must verify certain \ $\Sigma_k$ \ and \ $\Pi_k$ \ sentences which arise as he is playing the game. To
verify a \
$\Sigma_k$
\ sentence of the form \ $\exists v\eta(v)$ \ where \ $\eta$ \ is \ $\Pi_{k-1}$, \ player $\mathbf{I}$ must put the
formula \ $\eta(\iota v\sigma(v))$ \ in \ $T^*$ \ for some \ $\Sigma_{k-1}$ \ formula \ $\sigma$. \ However, such a
requirement (as stated) would present an open condition on player $\mathbf{I}$. Thus, to keep the requirements closed,
player $\mathbf{I}$ will be forced to bound the move at which he will verify a \ $\Sigma_k$ \ sentence in advance. This
is the purpose of player $\mathbf{I}$'s entry \ $m_j$, \ presented at move \ $j$. \ The integer \ $m_j$ \ is
player $\mathbf{I}$'s prediction of the move at which he will verify certain \ $\Sigma_k$ \ sentences.

Recall that \ $\theta_{\textup{sk}}(x,u,v,w_1,F)$ \ is the \ $\Sigma_{k-1}$ \ formula, in the parameters \ $w_1$ \ and
\ $F$, \ defining a \ $\Sigma_{k-1}$ \ Skolem function for \ $\overcore$. \ In addition, remember that \
$F$ \ is the \ $\le_{\BK}$--least such parameter.
\be
\item[R-6:] 
\be
\item The sentence 
\ $\forall v\forall w\forall x\forall
u(\theta_{\textup{sk}}(x,u,v,\ux_1,\uF)\land\theta_{\textup{sk}}(x,u,w,\ux_1,\uF)\rightarrow v=w)$ \ 
is in \ $T_0$.
\item If \ $\psi(v_0,v_1)$ \ and \ $\tau(u)$ \ are \ $\Sigma_{k-1}$ \ and the formula \ $\exists v\psi(v,\iota
u\tau(u))$ \ is in \ $T_j$, \ then the sentence
\ $\exists v(\psi(v,\iota u\tau(u))\land \theta_{\textup{sk}}(\ux_i,\iota u\tau(u),v,\ux_1,\uF))$ \ 
is in \ $T_{j+1}$ \ for some \ $i\in\omega$.
\item If the sentence \ $(\iota v\sigma(v)<_\BK\uF)$ \ is in \ $T_j$, \ then either the sentence
\[\exists v\exists w\exists x\exists
u(\theta_{\textup{sk}}(x,u,v,\ux_1,\iota v\sigma(v))\land\theta_{\textup{sk}}(x,u,w,\ux_1,\iota v\sigma(v))\land v\ne
w)\] is in \ $T_{m_j}$, \ or the sentence
\[\exists v\psi(v,\iota u\tau(u)) \land \lnot\exists x\exists v(\psi(v,\iota u\tau(u))\land
\theta_{\textup{sk}}(x,\iota u\tau(u),v,\ux_1,\iota v\sigma(v)))\] is in \ $T_{m_j}$ \ for some \ $\Sigma_{k-1}$ \ formulae
\ $\psi(v_0,v_1)$ \ and \ $\tau(u)$. 
\ee
\ee
\begin{rmk} Requirements R-6(a) and R-6(b) force player $\mathbf{I}$ to assert that \ $\theta_{\textup{sk}}$ \
with parameters \ $\uF$ \ and \ $\ux_1$ \ defines a \ $\Sigma_{k-1}$ \ Skolem function, whereas requirement R6-(c)  
compels player $\mathbf{I}$ to show that anything  $<_{\BK}$--less than \ $\uF$ \ does not define such a Skolem function.
\end{rmk}
Our next requirement will fix the interpretation of the constant symbol \ $\uG$. \ Recall that \ $\seq{\G,w_0}$ \
witnesses the fact that \ $\overcore$ \ realizes the type \ $\Gamma$. \ Also remember the
\ $\Sigma_k$ \ formula \ $\theta_i$ \ in the type \ $\Gamma$ \ used to define \ $Q_{i+1}$, \ and  the order \ $<_\Gamma$ \ 
on \ $\Gamma$ \ of order type \ $\omega$. \ Recall, as well, that in R-1 the constant \ $\ux_2$ \ is to be interpreted as the
real \ $w_0$.
\be
\item[R-7:] 
\be
\item If \ $\forall u \psi(u,v)$ \ is in the type \ $\Gamma$ \ where \ $\psi$ \ is \ $\Sigma_{k-1}$ \ and if \
$\tau(u)$ \ is any \ $\Sigma_{k-1}$ \ formula, then  the sentence \ $\lnot\psi(\iota v\tau(u), \seq{\uG,\ux_2})$ \ 
is not in \ $T_j$.
\item If \ $\exists u \psi(u,v)$ \ is in the type \ $\Gamma$ \ where \ $\psi$ \ is \ $\Pi_{k-1}$ \ and if \
$\exists u \psi(u,v)<_{\Gamma}\theta_j$, \ then for some \ $\Sigma_{k-1}$ \ formula \ $\tau(u)$ \ the sentence
\ $\psi(\iota v\tau(u), \seq{\uG,\ux_2})$ \ is  in \ $T_j$.
\item If the sentence \ $(\iota v\sigma(v)<_{\BK}\uG)$ \ is in \ $T_j$ \ for some \ $\sigma(v)$ \ in \
$\Sigma_{k-1}$, \ then either 
\be
\item there is a \ $\Sigma_{k-1}$ \ formula \ $\psi(u,v)$ \ such that 
the formula \ $\forall u\psi(u,v)$ \ is in the type \ $\Gamma$, \ but for some \ $\tau(u)$ \ in \ $\Sigma_{k-1}$ \
the sentence \ $\lnot\psi(\iota v\tau(u), \seq{\iota v\sigma(v),\ux_2})$ \ 
is in \ $T_{m_j}$, \ or 
\item there is a \ $\Pi_{k-1}$ \ formula \ $\psi(u,v)$ \ such that the formula \ $\exists u\psi(u,v)$ \ is one of the
first \ $m_j$ \ elements of \ $\Gamma$ \ under \ $<_{\Gamma}$, \ and the sentence \ $\psi(\iota v\tau(u), \seq{\iota
v\sigma(v),\ux_2})$
\  is not in \ $T^*$ \ for all \ $\Sigma_{k-1}$ \ formulae \ $\tau(u)$.
\ee
\ee
\ee
\begin{rmk} Requirements R-7(a) and R-7(b) oblige player $\mathbf{I}$ to produce a model in which the term \
$\seq{\uG,\ux_2}$ \ realizes the type \ $\Gamma$. \ Requirement R-7(c) forces player $\mathbf{I}$ to verify that the
constant symbol \ $\uG$ \ is the $<_{\BK}$--least such term.
\end{rmk}
Next, we shall fix the meaning of the constant symbols \ $\urho$ \ and \ $\up$. \ Recall the uniformly $\Sigma_1$ Skolem
function \ $h$ \ (see Definition~{\skolemdef} of \cite{Part1}).  

\be
\item[R-8:] 
\be
\item The sentence \ $(\urho\le\underline\kappa)$ \ is in \ $T_0$.
\item If \ $\tau(u)$ \ is \ $\Sigma_{k-1}$ \ and the sentence \ $\exists v(v=\iota u\tau(u))$ \ is in \ $T_j$, \ then
the sentence \ $(\iota v\sigma(v)\in\urho \land h(\ux_i,\seq{\iota v\sigma(v),\up})=\iota u\tau(u))$ \ is in \
$T_{j+1}$ \ for some \ $\Sigma_{k-1}$ \ formula \ $\sigma(v)$ \ and some \ $i\in\omega$.
\item If \ $\tau(v)$ \ and \ $\varphi(v)$ \ are \ $\Sigma_{k-1}$ \ and the sentence \ $(\iota v\tau(v)\in\urho\land
\iota v\varphi(v)\in [\OR]^{<\omega})$ \ is in \ $T_j$, \ then the sentence \
$\lnot(\exists x\in\uR)(\exists\gamma\in\iota v\tau(v))(h(x,\seq{\gamma,\iota v\varphi(v)})=\iota u\sigma(u))$
\ is in \ $T_{m_j}$ \ for some \ $\sigma(u)$ \ in \ $\Sigma_{k-1}$.
\item If \ $\tau(v)$ \ is \ $\Sigma_{k-1}$ \ and the sentence \ $(\iota v\tau(v)<_{\BK}\up)$ \ is in \ $T_j$, \ then
the sentence \ $\lnot(\exists x\in\uR)(\exists\gamma\in\urho)(h(x,\seq{\gamma,\iota v\tau(v)})=\iota
u\sigma(u))$ \ is in \ $T_{j+1}$ \ for some \ $\Sigma_{k-1}$ \ formula \ $\sigma(u)$.
\ee
\ee
\begin{rmk} Requirement R-8(b) coerces player $\mathbf{I}$ to prove that his model will satisfy the property that \
$h$, \ with parameter \ $\up$ \ and domain $\subseteq\uR\times\urho$, \ is  onto.   Requirement R-8(c) forces player
$\mathbf{I}$ to verify that the constant symbol \ $\urho$ \ is the smallest term that meets requirement R-8(b).
Similarly, requirement R-8(d) forces player $\mathbf{I}$ to verify that the constant symbol \ $\up$ \ is the $<_{\BK}$
\ smallest term that meets requirement R-8(b). Definition~{\sounddef} and Lemma~{\soundlemma} of \cite{Part1}
imply that these requirements will provide an accurate interpretation of the constant symbols \ $\urho$ \ and \
$\up$.
\end{rmk}
Our final requirement will make player $\mathbf{I}$ correctly interpret the constant symbols \ $\uQ_i$. \ Recall the
sequence \ $\seq{\theta_i : i\in\omega}$ \ of \ $\Sigma_k$ \ formula in
\ $\Gamma$ \ used to define the \ $Q_i$'s. \ For each such \ $\theta_i(v)$ \ let \ $\psi_i(u,v)$ \ be the \ $\Pi_{k-1}$
\ formula so that \ $\theta_i(v)=\exists u\psi_i(u,v)$.
\be
\item[R-9:] 
\be
\item The sentence \ $(\uQ_0=\seq{\urho,\up})$ \ is in \ $T_0$.
\item For any integer \ $i>0$ \ there is \ $\Sigma_{k-1}$ \ formula \ $\tau(v)$ \ with support \ $i$ \ such that the
sentence \ $\psi_{i-1}(\iota v\tau(v),\seq{\uG,\ux_2})$ \ is in \ $T_i$.
\item For all \ $i>0$, \ all \ $j\in\omega$ \ and all \ $\sigma(v)$ \ in \ $\Sigma_{k-1}$, \ if the sentence \ $(\iota
v\sigma(v)<_{\BK}\uQ_i)$ \ is in \ $T_j$, \ then the sentence \ $\psi_{i-1}(\iota u\tau(u,\iota
v\sigma(v)),\seq{\uG,\ux_2})$
\ is not in \ $T^*$ \ for all \ $\Sigma_{k-1}$ \ formulae \ $\tau(u,v)$ \ with support \ $i-1$.
\ee
\ee
\begin{rmk} Requirements R-9(b) and R-9(c) require player $\mathbf{I}$ to verify that the definition of the sequence \
$\seq{\theta_i : i\in\omega}$ \ is correct and also ensure that player $\mathbf{I}$ accurately interprets the elements of
the sequence \ $\seq{Q_i : i\in\omega}$.
\end{rmk}
We shall now define when an initial position \ $u$ \ of the game \ $G_x$ \ is honest. Let 
\[u=\seq{\seq{T_j,s_j,\f_j,m_j,y_j}:j<\n}\]
be a position of length \ $\n$, \ where \ $\n\in\omega$, \ and let \ $\m\in\omega$ \ be such that
\[\seq{x_i:i<\m}=s_0{^\frown}\seq{y_0}{^\frown}s_1{^\frown}\seq{y_1}{^\frown}\cdots
{^\frown}s_{\n-1}{^\frown}\seq{y_{\n-1}},\]
and define \ $I_u$, \ an initial interpretation in \ $\overcore$ \ of the constant symbols, as follows:
\[\begin{array}{l}
I_u(\ux_i)=x_i \text{ \ for $i<\m$},\\
I_u(\uF)=F, \ I_u(\uG)=\G,\\
I_u(\urho)=\rho_{\overcore}^{}, \ I_u(\up)=p_{\overcore}^{},\\ 
I_u(\uQ_i)=Q_i \text{ \ for all $i\in\omega$}.\\
\end{array}\]
We shall call the initial position \ $u$ \ {\it reasonable\/} if it is not an immediate loss for player $\mathbf{I}$ 
because of \ $\seq{T_j : j<\n}$ \ in the following sense: all of the conditions on \ $T_j$ \ in requirements  R-1
through R-4 and requirements R-6 through R-9 are satisfied, for each \ $j<\n$. \ We can now characterize the 
winning positions for player $\mathbf{I}$ in \ $G_x$ \ in which he has been honest about the model he is constructing.
Recall that \ $\theta(v)$ \ is the \ $\Sigma_k$ \ formula in the language \ $\Lng_\overn$, \ containing only the
parameters \ $\y, F, \G, Q_0, Q_1, \dots, Q_\ell$, \ that defines \ $P$ \ in \ $\H^{i_0}$.

We say that the initial position \ $u$, \ of length \ $n$, \ is  {\it $x$--honest\/} if the following eight conditions are
satisfied:
\be
\item[H-1:] $\H^{i_0}\models \theta(x)$.
\item[H-2:] The initial position \ $u$ \ is reasonable.
\item[H-3:] If \ $\n>0$, \ then \ $s_0(0)=x$, \ $s_0(1)=w_1$, \ $s_0(2)=w_0$, \ $s_0(3)=\y$.
\item[H-4:] $(\overcore, I_u)\models\bigcup\limits_{j<\n} T_j$.
\item[H-5:] Let \ $\sigma_0,\dots, \sigma_i, \dots, \sigma_m$ \ enumerate those \ $\B_{k-1}$ \ sentences that satisfy
the following three properties:  
\be
\item[(1)] $n(\sigma_i)<\n$
\item[(2)] $\sigma_i$ \ has the form \ $(\iota v\phi_i(v)\in\F_{n_i})$ \ where \ $\phi_i$ \ is \ $\Sigma_{k-1}$ \ and \
$n_i\in\omega$ 
\item[(3)] $(\overcore,I_u)\models\sigma_i$. 
\ee
For each \ $i\le m$ \ let \ $\h_i\in\F^{\overcore}_{n_i}$ \ be such that \
$(\overcore, I_u)\models\phi_i(\h_i)$. \ Let \ $\tau=\{(\h_i,\f_{n(\sigma_i)}): i\le m\}$. \ Then there is \ $\Phi\colon
\F^{\overcore}  \leadsto \F^{\overcore}$ \ such that \ $\tau\subseteq\Phi$.
\ee
\begin{rmk} For a definition of \ $\Phi\colon \F^{\overcore}  \leadsto \F^{\overcore}$, in condition H-5,  \ 
see Definition~{\quasi} of \cite{Part1}.
\end{rmk}
The following three conditions for  $x$--honesty guarantee that player $\mathbf{I}$ has made predictions \ $m_j$ \
which he can fulfill (see requirements R-6 through R-8).
\be
\item[H-6:] If \ $j<\n$ \ and the sentence \ $(\iota v\sigma(v)<_{\BK}\uF)$ \ is in \ $T_j$ \ for some \
$\Sigma_{k-1}$ \ formula \ $\sigma$, \ then either 
\be
\item[(i)] $(\H^{m_j},I_u)\models\exists v\exists w\exists x\exists
u(\theta_{\textup{sk}}(x,u,v,\ux_1,\iota v\sigma(v))\land\theta_{\textup{sk}}(x,u,w,\ux_1,\iota v\sigma(v))\land v\ne
w)$, \ or 
\item[(ii)] $(\H^{m_j},I_u)\models\exists u(\exists v\psi(v,u) \land \lnot\exists x\exists
v(\psi(v,u)\land
\theta_{\textup{sk}}(x,u,v,\ux_1,\iota v\sigma(v))))$ \ for some \ $\Sigma_{k-1}$ \ formula \ $\psi(v,u)$.
\ee
\item[H-7:] If \ $j<\n$ \ and the sentence \ $(\iota v\sigma(v)<_{\BK}\uG)$ \ is in \ $T_j$ \ for some \
$\Sigma_{k-1}$ \ formula \ $\sigma$, \ then either 
\be
\item[(i)] $(\H^{m_j},I_u)\models\lnot\varphi(\seq{\iota v\sigma(v),\ux_2})$ \ for some \ $\Pi_k$ \ formula \ $\varphi$ \ in
the type \ $\Gamma$, \ or 
\item[(ii)] $(\H^{m_j},I_u)\models\lnot\varphi(\seq{\iota v\sigma(v),\ux_2})$ \ for some \ $\Sigma_k$ \ formula \ $\varphi$ \
which is one of the first \ $m_j$ \ elements in the type \ $\Gamma$.
\ee
\item[H-8:] 
\be
\item If \ $j<\n$ \ and the sentences \ $(\iota v\tau(v)\in\urho)$ \ and \ $(\iota v\varphi(v)\in[\OR]^{<\omega})$ \
are in \ $T_j$ \ for some \ $\Sigma_{k-1}$ \ formulae \ $\tau$ \ and \ $\varphi$, \ then
\[(\H^{m_j},I_u)\models\lnot(\exists x\in\uR)(\exists\gamma\in\iota v\tau(v))(h(x,\seq{\gamma,\iota
v\varphi(v)})=\iota u\sigma(u))\] 
for some \ $\sigma(u)$ \ in \ $\Sigma_{k-1}$.
\item If \ $j<\n$ \ and the sentence \ $(\iota v\tau(v)<_{\BK}\up)$ \ is in \ $T_j$ \ for some \ $\Sigma_{k-1}$ \
formulae \ $\tau$, \ then
\[(\H^{m_j},I_u)\models\lnot(\exists x\in\uR)(\exists\gamma\in\urho)(h(x,\seq{\gamma,\iota v\tau(v)})=\iota
u\sigma(u))\] \ for some \ $\Sigma_{k-1}$ \ formula \ $\sigma(u)$.
\ee
\ee
This completes our description of \ $x$--honesty. The assumption \ $\overcore\models\AD$
\ will be used explicitly in the proof of the following claim.
\setcounter{clm}{0}
\begin{clm}\label{claim1} The set \ $\{(x,u) : u \text{ is an $x$--honest position of length  $\n$}\}$ \ is \
$\boldface{\Sigma}{\omega}(\H^{\max\{i_0,\n\}})$.
\end{clm}
\begin{proof}[Proof of Claim~\ref{claim1}] Conditions H-1, H-2, and H-3 are clearly \
$\boldface{\Sigma}{\omega}(\H^{\max\{i_0,\n\}})$. \ Condition H-4 is first order over \ $\H^{\max\{i_0,\n\}}$ \ because
of our restrictions on the sentences in \ $T_j$ \ for \ $j<\n$. \ Since  \  $\H^{\max\{i_0,\n\}}\prec_1\overcore$ \ and  \
$\tau$ \ is an element of the structure \ $\H^{\max\{i_0,\n\}}$, \ Theorem~{\defontau} of \cite{Part1} implies that 
condition H-5 is \ $\boldface{\Sigma}{\omega}(\H^{\max\{i_0,\n\}})$. \ Since the proofs dealing with conditions H-6 and
H-8 are simpler than the argument addressing condition H-7, we shall just
prove that H-7 is \ $\boldface{\Sigma}{\omega}(\H^{\max\{i_0,\n\}})$. \ Consider the relation \ $R(K,\varphi)$ \ defined by
\[R(K,\varphi)\iff(\exists m\ge n)\left[\H^n\models(K<_{\BK}\G) \land \varphi\in\Gamma \land \H^m\models
\varphi(\seq{K,w_0})\right].\]
If we can show that \ $R$ \ is \ $\boldface{\Sigma}{\omega}(\H^n)$, \ then it
is straightforward to verify that condition H-7 is first order over \ $\H^{\max\{i_0,\n\}}$. \ 
Recall that for each \ $m\in\omega$  \ we have that \ $\H^m\prec_1\overcore$. \ Therefore, the satisfaction relation \
$\H^m\models\varphi(\seq{K,w_0})$ \ on \ $\Sigma_k$ \ formula \ $\varphi$ \ with parameter \ $K$, \ is uniformly \
$\Sigma_\omega(\H^m)$ \ (we are `equating' each \ $\varphi$ \ with a G\"odel number). It follows that \ $R\subseteq
H^n$ \ is in \ $\boldface{\Sigma}{0}^\omega(\overcore,H^0,H^1,\dots)$. \ Since \ $\overcore\models\AD$, \ Theorem
\ref{ADdefinable} implies that \ $R$ \ is in \ $\boldface{\Sigma}{\omega}(\H^{n})$ \ (in this case, \ $e(n)=n$). This completes
the proof of Claim 1.
\end{proof}
\begin{clm}\label{claim2} For all \ $x\in\R$ \ and all \ $u$, \ the following are equivalent:
\be
\item $u$ \ is $x$--honest
\item $u$ \ is a winning position for player $\mathbf{I}$ in $G_x$.
\ee
\end{clm}
\begin{proof}[Proof of Claim~\ref{claim2}] We shall first prove that \ $(1)\Rightarrow(2)$ \ and then show that \
$(2)\Rightarrow(1)$. 

$(1)\Rightarrow(2)$: \ Let \ $u$ \ be an  $x$--honest position of length \ $\n$. \ We shall show that
\begin{equation}(\exists T,s,\f,m)(\forall y\in\R)( u^{\frown}\seq{T,s,\f,m,y}\text{ is
$x$--honest}).\tag{$\star$}\end{equation} Since \ $u$ \ is an arbitrary $x$--honest position, $(\star)$ implies that
player $\mathbf{I}$ can win the game \ $G_x$ \ by repeatedly playing honest positions. To establish $(\star)$, we shall
assume that player $\mathbf{I}$ and his opponent have produced \ $u$ \ while playing the game \ $G_x$. \ It is now player
$\mathbf{I}$'s move and we shall show that he can continue to play honestly. We note that since
\ $u$ \ is $x$--honest, we have that \ $u$ \ satisfies conditions H-2 and H-4; thus, player $\mathbf{I}$ has expressed as much
of the truth of \ $\overcore$ \ as he was required to tell. Because \ $u$ \ meets conditions H-6, H-7 and H-8, player
$\mathbf{I}$ has made predictions \ $m_j$, \ for each \ $j<n$, \ which he can fulfill. Therefore, player $\mathbf{I}$ can
choose \ $T$ \ and \ $s$ \ so that the new position will satisfy conditions H-2 and H-4 (if \ $n=i_0$, \ then one also needs
the fact that \ $u$ \ satisfies condition H-1). Continued satisfaction of H-1 and H-3 is easy to assure. By choosing \ $m$ \
large enough, conditions H-6, H-7 and H-8 can be fulfilled. Finally player $\mathbf{I}$ must choose \ $\f$ \ so as to make
certain that condition H-5 is satisfied. Since \ $u$ \ is $x$--honest, let \ $\tau$ \ and \ $\Phi\colon
\F^{\overcore}  \leadsto \F^{\overcore}$ \  be as stated in condition H-5. Because \ $u$ \ is reasonable and because of the
restrictions in H-5, it follows that \ $\tau$ \ is in \ $\H^n$. \ Hence, \ $\tau\subseteq\F^{\H^{n}}\times\F^{\H^{n}}$ \ and
\ $\tau\subseteq\Phi$. \ Since \ $\H^{n+1}\prec_1\overcore$, \ Theorem 3.25 of \cite{Cfsrm} implies that is a \ $\Phi'\colon
\F^{\H^{n+1}}  \leadsto \F^{\H^{n+1}}$ \ where \ $\tau\subseteq\Phi'$. \ Thus, for any \ $\h\in\F^{\H^{n+1}}$ \ there is an \
$\f\in\F^{\H^{n+1}}$ \ so that \ $\tau'=\tau\cup\{\seq{\h,\f}\}\subseteq\Phi'$. \ Theorem 3.25 of \cite{Cfsrm} then implies
that there is a \ $\Phi^*\colon\F^{\overcore}  \leadsto \F^{\overcore}$ \ such that \ $\tau'\subseteq\Phi^*$. \ So,
continued satisfaction of condition H-5 can be guaranteed. Therefore, $(\star)$ has been established.

$(2)\Rightarrow(1)$: \ Suppose that \ $\Sigma$ \ is a winning strategy for player $\mathbf{I}$ starting from \ $u$. \ We shall
prove that \ $u$ \ is $x$-honest by using a `generic run' argument. This technique was used in the proof of
Lemma 4.7 of \cite{Crcm}. So, let \ $H=\seq{\seq{T_j,s_j,\f_j,m_j,y_j}:j<\omega}$ \ be a generic run, according to \
$\Sigma$, \ in the game \ $G_x$. \ Let \
\ $\seq{x_i:i\in\omega}=s_0{^\frown}\seq{y_0}{^\frown}s_1{^\frown}\seq{y_1}{^\frown}\cdots$ \ 
and let \ $T^*=\bigcup\limits_{j\in\omega}T_j$. \ By requirement R-2(a), the theory \ $T\cup T^*$ \ is consistent. Let \
$\fB$ \ be a model of \ $T\cup T^*$. \ By payoff requirements R-6(a) and R-6(b) we have that \
$\fA\prec_{k-1}\fB$ 
\ where 
\[\fA=\Hull_{k-1}^{\fB}(\{\ux_i^\fB:i\in\omega\}\cup\{\uR^\fB,\ul{\kappa}^\fB,\uF^\fB, \uG^\fB, \urho^\fB,
\up^\fB\}\cup\{\uQ_i^\fB:i\in\omega\}).\]
By requirement R-2(d), we have that \ $\fA\models\T^{\overn}$ \ (recall that \ $\overn=n(\core)$). 
Requirement R-4 implies that \ $\uR^\fA=\{\ux_i^\fB:i\in\omega \}$. \ For \ $\h\in\F_n^\fA$, \ let \ $\tau$ \ be the first
formula of the form \ $(\iota v\phi(v)\in\F_n)$, \ in the enumeration induced by the generic run \ $H$, \ such that 
\ $B\models(\h=\iota v\phi(v)\land \h\in\F_n)$ \ 
and define \ $\sigma(\h)=\f_{n(\tau)}$. \ Requirement R-5 implies that the function \ $\sigma\colon\F^{\fA}\to\F^{\overcore}$ \
is  $\ul{E}_\overn$--extendible.   So, by Theorem~{\criterion} of \cite{Part1}, there is an $\overn$--iterable
premouse \ $\kouse$ \ such that \ $\kouse^\overn$ \ is the transitive collapse of \ $\fA$. \ Hence, \
$\fA$ \ is isomorphic to \ $(\kouse^\overn,I)$ \ for some interpretation \ $I$ \ of the constants. By requirement R-8, \
$\kouse$ \ is critical and \ $\rho_{\kouse}^{\overn+1}\le\kappa^\kouse$. \ Therefore, \ $\kouse$ \ is a mouse with \
$n(\kouse)=\overn=n(\core)$. \ We shall write \ $\overkouse=\kouse^\overn$. \ Now, by genericity we have that \
$\R^{\overkouse}=\R=\R^V$ \ where \ $V$ \ is the ground model over which \ $H$ \ is generic. Requirements R-7(a)
and R-7(b) imply that \ $\overkouse$ \ realizes the type \ $\Gamma$. \ Hence, \ $\overkouse$ \ realizes the type \
$\overline{\Upsilon}^*$. \ Because \ $\sigma\colon\F^{\overkouse}\to\F^{\overcore}$ \ is  $\ul{E}_\overn$--extendible, Lemma
\ref{typecomparison2} implies that \ $\overkouse=\overcore$. \ Finally, payoff requirements R-2(e), R-6, R-7, R-8 and R-9
ensure that 
\[\begin{array}{l}
I(\ux_i)=x_i \text{ \ for all $i\in\omega$},\\
I(\uF)=F, \ I(\uG)=\G,\\
I(\urho)=\rho_{\overcore}^{}, \ I(\up)=p_{\overcore}^{},\\ 
I(\uQ_i)=Q_i \text{ \ for all $i\in\omega$}.\\
\end{array}\]
It is now straightforward to verify that the conditions of $x$--honesty hold for \ $u$ \ in the extension \ $V[H]$. \
All of the conditions in the definition of $x$-honesty, except condition H-5, are easily shown to be absolute between \ $V$ \
and \ $V[H]$. \ To show that H-5 is absolute,  we note that the map \ $\sigma\colon\F^{\overcore}\to\F^{\overcore}$,
\ after identifying \ $\fA$ \ and \ $\overcore$, \ is in \ $V[H]$ \ (but not in \ $V$)  and can be used in \ $V[H]$ \ to
verify that condition H-5 holds for \ $u$. \ Thus, there is a \ $\Phi\colon \F^{\overcore}  \leadsto \F^{\overcore}$ \ in \
$V[H]$ \ that verifies condition  H-5 for \ $u$. \ However, Lemma 3.19 of \cite{Cfsrm} and its proof imply that, if there
exists such an \ $\Eq$--extendible quasi-map  extending \ $\tau$, \ then there is a \ $\boldface{\Sigma}{\omega}(\overcore)$
\ such quasi-map (see \cite[Definition 3.20]{Cfsrm}). Therefore, condition H-5 is absolute between \ $V$ \ and \ $V[H]$. \ 
It follows that \ $u$ \ is $x$--honest in \ $V$ \ and this completes the proof of Claim~\ref{claim2}.
\end{proof} 
Claim~\ref{claim2} applied to the empty position implies that \ $x\mapsto G_x$ \ is a closed game representation of \ $P$. \
Claims~\ref{claim1} and \ref{claim2} \ imply, as stated at the beginning of the proof of Lemma~\ref{closedk>1}, that the
resulting Moschovakis scale  on \ $P$ \ is \ $\boldface{\Sigma}{k}(\overcore)$. \ This
completes the proof of Lemma~\ref{closedk>1}.
\end{proof}
Lemma~\ref{closedk>1} asserts that if the set of reals \ $P$ \ is \ $\Sigma_k(\H^{i_0})$ \ in the parameters $\y,F, \G, Q_0,
\dots, Q_\ell$, \ then \ $P$ \ has a \ $\boldface{\Sigma}{k}(\overcore)$ \ scale. We can now prove this will hold when 
one allows arbitrary parameters.
\begin{lemma} If \ $\overcore\models\AD$ \ and \ $P\subseteq\R$ \ is \ $\boldface{\Sigma}{k}(\H^{i})$ \ for some \
$i\in\omega$, \ then \ $P$ \ has a scale which is \ $\boldface{\Sigma}{k}(\overcore)$.
\end{lemma}
\begin{proof} Suppose that \ $P\subseteq\R$ \ is \ $\boldface{\Sigma}{k}(\H^{i})$ \ for some  \ $i\in\omega$. \ From the
definition of \ $\H^i$ \ and because \ $\H^i$ \ has a \ $\Sigma_{k-1}$ \ Skolem function in the parameter \ $F$, \ it follows
that for some real \ $y$ \ there is a  \ $\Sigma_k$ \ formula in the parameters \ $\y,F, \G, Q_0, Q_1, \dots, Q_i$ \ that
defines \ $P$ \ in \ $\H^{i}$. \ Lemma~\ref{closedk>1} implies there is a scale on \ $P$ \ which is
\ $\boldface{\Sigma}{k}(\overcore)$.
\end{proof}
This completes Case 1 in our proof of Theorem~\ref{goodcovering}.

{\bf Case 2: $\mathbf{k=1}$ and $\mathbf{\widehat{\OR}^{\overcore}}$ is a limit ordinal.} \ Using \ $\G$, \ $w_0$ \ and
an argument similar to the proof of Lemma~\ref{hull}, we can now define a canonical sequence \ $\seq{\H^i :
i\in\omega}$ \ of \ $\Sigma_{\omega}$ \ hulls. At the same time we shall also define a sequence \ $\seq{\theta_i :
i\in\omega}$ \ of \ $\Sigma_1$ \ formula in \ $\Gamma$, \ and a sequence \ $\seq{\beta_i : i\in\omega}$ \ of ordinals \
$\beta_i\in\widehat{\OR}^{\overcore}$. \ Recall, \ $\seq{\G,w_0 }$ \ is the witness verifying that \ $\overcore$ \
realizes the type \ $\Gamma$. \ Since \ $k=1$, \ $\Gamma$ \ consists of \ $\Sigma_1$ \ and \ $\Pi_1$ \ formulae.

First, define \ $\beta_0$ \ be the least ordinal \ $\beta$ \ such that \
$\kappa^{\overcore}<\beta<\widehat{\OR}^{\overcore}$ \ and \
$p^{}_{\overcore}\in\overcore^{\,\beta}$. \ Note that \ $\rho^{}_{\overcore}\in\overcore^{\,\beta_0}$. \ Assuming that
$\beta_0, \beta_1, \dots, \beta_i$ \  have been defined, we define \ $\theta_i$, \ $\beta_{i+1}$ \ and  $\H^i$ \ as follows:

\smallskip
\noindent{\sl Definition of \ $\theta_i$\/}: \ Since \ $\overcore^{\beta_i}$ \ does not realize \ $\Gamma$,
\ let \ $\theta_i$ \ to be the \ $<_\Gamma$\,--least \ $\Sigma_1$ \ formula \ $\theta(v)\in\Gamma$ \ such
that \ $\overcore^{\,\beta_i}\not\models\theta(\seq{\G,w_0})$.

\smallskip
\noindent{\sl Definition of \ $\beta_{i+1}$\/}: \ Because \ $\theta_i$ \ is \ $\Sigma_1$, \
there is an ordinal \
$\beta\in\widehat{\OR}^{\overcore}$ \ such that
\begin{equation}
\overcore^{\,\beta}\models\theta_i(\seq{\G,w_0}).
\label{makeit}\end{equation} 
Let \ $\beta_{i+1}$ \ be the least such \ $\beta>\beta_i$ \ satisfying (\ref{makeit}).

\smallskip
\noindent{\sl Definition of \ $\H^i$\/}: \ Define
\ $\H^i=\Hull_{\omega}^{\overcore^{\,\beta_{i+1}}}(\R\cup\{\G,
\rho^{}_{\overcore},p^{}_{\overcore},\beta_0, \beta_1, \dots,
\beta_i\})$. 

\begin{lemma}\label{lemma_k=1limit} Let \ $\seq{\H^i :
i\in\omega}$, \ $\seq{\theta_i :
i\in\omega}$ \ and \ $\seq{\beta_i :
i\in\omega}$ \ be as above. Then \ $\bigcup\limits_{i\in\omega}\H^i=\overcore$ \ and for
all \ $i\in\omega$
\be
\item[(a)] $\H^i\prec_{0}\H^{i+1}\prec_{0}\overcore$
\item[(b)] $\mathcal{H}^i$ \ is in \ $\boldface{\Sigma}{1}(\mathcal{H}^{i+1})$
\item[(c)] $\pow(\R)\cap\boldface{\Sigma}{\omega}(\H^i)\subseteq\overline{C}$
\item[(d)] there exists a \ $\boldface{\Sigma}{1}(\mathcal{H}^{i+1})$ \ function \ $f\colon\R\maps{onto}H^i$
\item[(e)] $\H^i\not\models\theta_i(\seq{\G,w_0})$ \ and \ $\H^{i+1}\models\theta_i(\seq{\G,w_0})$.
\ee
\end{lemma}
\begin{proof} For each \ $i\in\omega$, \ items (a)--(e) are clear; for example, since \
$\overcore^{\,\beta_{i}}\prec_0\overcore^{\,\beta_{i+1}}$, \ it follows that \ $\H^i\prec_{0}\H^{i+1}$. \ To prove that \
$\bigcup\limits_{i\in\omega}\H^i=\overcore$, \ let \ $\H=\bigcup\limits_{i\in\omega}\H^i$. \ Because the sequence \
$\seq{\beta_i : i\in\omega}$  \ is cofinal in \ $\widehat{\OR}^{\overcore}$, \ it follows that 
\ $\H\prec_1\overcore$. \ Arguing as in the proof of Lemma 2.29 of \cite{Cfsrm}, there is a mouse \ $\kouse$ \ such
that \ $n(\kouse)=n(\core)=\overn$ \ and there is an  isomorphism \ $\pi\colon\H\to\overkouse$ \ where \ $\pi$ \
is the transitive collapse of \ $\H$. \ Therefore, \ $\overkouse$ \ realizes the type \ $\overline{\Upsilon}^*$. \
Since \ $\pi^{-1}\colon\overkouse\mapsigma{1}\overcore$, \ Lemma~\ref{typecomparison1} implies that \ $\overkouse=\overcore$.
\ Thus, \ $\pi\colon\H\to\overcore$ \ is an (onto) isomorphism. The following two claims establish that \ $\H=\overcore$.
\setcounter{clm}{0}
\begin{clm} Assume that \ $\pi(\G)=\G$, \ $\pi(\rho^{}_{\overcore})=\rho^{}_{\overcore}$, \
$\pi(p^{}_{\overcore})=p^{}_{\overcore}$ \ and \
$\pi(\beta_i)=\beta_i$ \ for all \ $i\in\omega$.
\ Then \
$\H=\overcore$. 
\end{clm}
\begin{proof}[Proof of Claim 1] Clearly, \ $\pi(x)=x$ \ for all \ $x\in\R$. \ Since \
\[\H=\bigcup\limits_{i\in\omega}\H^i=\Hull_{1}^{\H}(\R\cup
\{\G,\rho^{}_{\overcore},p^{}_{\overcore},\beta_0,\dots,\beta_i,\dots\}),\] it follows that 
\[\overcore=\Hull_{1}^{\overcore}(\R\cup
\{\pi(\G),\pi(\rho^{}_{\overcore}),\pi(p^{}_{\overcore}), \pi(\beta_0),\dots,\pi(\beta_i),\dots\}).\]
Therefore, our assumption implies that \ $\H=\overcore$.
\end{proof}

\begin{clm} $\pi(\G)=\G$, \ $\pi(\rho^{}_{\overcore})=\rho^{}_{\overcore}$, \
$\pi(p^{}_{\overcore})=p^{}_{\overcore}$ \ and \
$\pi(\beta_i)=\beta_i$ \ for all \ $i\in\omega$.
\end{clm}
\begin{proof}[Proof of Claim 2] The proof of this claim is, for the most part, a repetition of the argument used to
establish Claim 2 in the proof of Lemma~\ref{lemma_k>1}.
\end{proof}
This completes the proof of the lemma. 
\end{proof}

Letting \ $e(i)=i+1$ \ for all \ $i\in\omega$, \ it then follows that  \ $\langle\H^i : i\in\omega\rangle$ \
is a good covering of \ $\overcore$.

\begin{lemma}\label{onto2} There is a partial \ $\boldface{\Sigma}{1}(\overcore)$ \ map of \ $\R$ \ onto \ $\overC$, \ the
domain of \ $\overcore$.
\end{lemma}
\begin{proof} Let \  $\H=\Hull_1^{\overcore}(\R\cup\{\G\})$. \ It follows that \ $\seq{\G,w_0}$ \ realizes the type \
$\overline{\Upsilon}^*$ \ in \ $\H$. \ Arguing as in the proof of Lemma 2.29 of \cite{Cfsrm}, there is a mouse \ $\kouse$ \
such that \ $n(\kouse)=n(\core)=\overn$ \ and there is an  isomorphism \
$\pi\colon\H\to\overkouse$ \ where \ $\pi$ \ is the transitive collapse of \ $\H$. \ Therefore, \ $\overkouse$ \ realizes
the type \ $\overline{\Upsilon}^*$. \ Since \ $\pi^{-1}\colon\overkouse\mapsigma{1}\overcore$, \ Lemma~\ref{typecomparison1}
implies that \ $\overkouse=\overcore$. \ Thus, \ $\pi\colon\H\to\overcore$ \ is an (onto) isomorphism. Because of the
definition of \ $\G$, \ we have that \ $\pi$ \ is the identity map. Hence, \ $\H=\overcore$ \ and we must have that the \
$\Sigma_1$ \ Skolem function \ $h_{\seq{\G,w_0}}$ \ (in the parameter \ $\seq{\G,w_0}$) \ maps \ $\R\times\R$ \ onto \ $C$. 
\end{proof}

\begin{lemma}\label{closedk=1limit} Consider the parameters \ $\y, \G, \rho^{}_{\overcore}, p^{}_{\overcore}, \beta_0, \beta_1,
\dots, \beta_\ell$ \ (as defined above) where \ $\y$ \ is a fixed real.  Let \ $\theta$ \ be a \ $\Sigma_1$ \ formula in the
language \ $\Lng_\overn$ \ augmented by these parameters. Fix \ $i_0\ge\ell$ \ and let \ $P\subseteq\R$ \ be defined by \
$P(x)$ \ if and only if \ $\H^{i_0}\models \theta(x)$. \ Then \ $P$ \ has a closed game representation and if \ $\overcore\models\AD$, \ then there
is a \ $\boldface{\Sigma}{1}(\overcore)$ \ scale on \ $P$.
\end{lemma}
\begin{proof}[Sketch of Proof] Since the main ideas of the proof of this lemma are the same as those used in the proof of
Lemma~\ref{closedk>1}, we shall only outline the relevant details needed to provide an explicit proof. 
 We want to construct a closed game representation of \ $P$ \ simple enough to yield a
\ $\boldface{\Sigma}{1}(\overcore)$ \ scale on \ $P$. \ Our closed game representation \ $x\mapsto G_x$ \
of \ $P$ \ will have the following property: If \ $P_n$ \ is defined by
\[P_n(x,u)\iff \textup{$u$ is a winning position for player $\mathbf{I}$ in $G_x$ of length $n$},\]
then \ $P_n$ \ is \ $\boldface{\Sigma}{\omega}(\H^{\max(i_0,n)})$. \ Assuming \ $\overcore\models\AD$, \ Lemma
\ref{poffinmouse} and Lemma~\ref{lemma_k=1limit}(c) will then imply that the Moschovakis scale \ $\seq{\varphi_i : i \in
\omega}$ \ on \ $P$ \ exists. Furthermore, the prewellordering \ $\le_i$ \ of \ $\R$ \ induced by \ $\varphi_i$ \ is also \
$\boldface{\Sigma}{\omega}(\H^{\max(i_0,i)})$. \ Lemma~\ref{lemma_k=1limit}(c) implies that \ $\le_i\,\in\overC$, the domain
of \ $\overcore$, \ for all \ $i\in\omega$. \ Lemma~\ref{onto2} gives a partial \
$\boldface{\Sigma}{1}(\overcore)$ \ map of \ $\R$ \ onto \ $\overC$. \ It follows that
any countable subset of \ $\overC$ \ is \ $\boldface{\Sigma}{1}(\overcore)$. \ Therefore, \ $\seq{\varphi_i : i \in
\omega}$ \ is a \ $\boldface{\Sigma}{1}(\overcore)$ \ scale on \ $P$. 

So, to prove the lemma, it suffices to construct a closed game representation \ $x\mapsto G_x$ \ so that \ $P_n$ \ is
first order definable over the structure \ $\H^{\max(i_0,n)}$, \ for each \ $n\in\omega$. \ In our construction of \ $G_x$ \ we want to 
force player $\mathbf{I}$ to describe the truth
in \ $\overcore$ \ about \ $\G$, $\rho^{}_{\overcore}$, $p^{}_{\overcore}$, and the \ $\beta_i$'s. \ To ensure that each \
$P_n$ \ is \ $\boldface{\Sigma}{\omega}(\H^{\max(i_0,n)})$ \ we shall impose some restrictions on player $\mathbf{I}$'s moves. In this case,
player $\mathbf{I}$ must describe a model in the language 
\[\Lng=\Lng_\overn\cup\{\uG,\urho,\up\}\cup\{\ux_i,\ul{\beta}_i,\ul{\overcore^{\beta_i}}: i\in\omega\}\]
where \ $\overn=n(\core)$. \ 
Player $\mathbf{I}$ must play a consistent and complete set of \ $\Sigma_{0}$ \ sentences in the language \ $\Lng$,  \
mentioning at move \ $j$ \ no sentences involving the constants \ $\ul{\beta}_i$ \ and \ $\ul{\overcore^{\beta_i}}$ \ for \
$i>j$. \ For each \ $j\in\omega$ \ player $\mathbf{I}$ must play at move \ $j$ \ the \ $\Sigma_0$ \ sentence \
$(\ul{\overcore^{\beta_j}}\models\PM)$, \ where we recall that the theory \ $\PM$ \ (see Definition~{\PMdefined} of
\cite{Part1}) can be axiomatized by a single sentence. At move \ $i_0+1$ \ player $\mathbf{I}$ must assert that some object,
definable over \ $\ul{\overcore^{\beta_{i_0+1}}}$ \ from the constants \ $\uG,\urho,\up,\ul{\beta}_{0},\dots,
\ul{\beta}_{i_0}$ \ and some ``real'' \ $\ux_k$ \ (which he has played), witnesses that the \ $\Sigma_1$ \ statement \
$\theta(\ux_0)$ \ will hold in the model he is constructing.  To prove that his model is mouse iterable, player $\mathbf{I}$
must play  `functions' from \ $\F^j$ \ where \ $\F^j=\F^{\H^j}$, \ at each move \ $j$. \  A typical run of the game has the
form
\begin{equation}
\begin{aligned}[c]
{}&{\mathbf{I}}\phantom{{\mathbf{I}}} \qquad T_0,s_0,\f_0,m_0 \qquad \phantom{y_0} \qquad T_1,s_1,\f_1,m_1 \qquad
\phantom{y_1}\quad
\\ {}&{\mathbf{I}}\mathbf{I} \qquad \phantom{T_0,s_0,\f_0,m_0}\qquad y_0 \qquad \phantom{T_1,s_1,\f_1,m_1} \qquad y_1 \quad 
\end{aligned}
\begin{gathered}[c]
{\cdots}
\end{gathered}
\tag*{$G_{x}$}\end{equation}
 where 
\bi
\item $T_j$ \ is a finite set of sentences each of which is in \ $\Sigma_0$ \ and has support \ $j$
\item $s_j\in\R^{<\omega}$, \ $\f_j\in \F^{j}$, \ $j<m_j\in\omega$, \ and \ $y_j\in\R$.
\ei
A \ $\Sigma_0$ \ formula \ $\varphi$ \ has {\it support\/} $n$ \ if it contains no constants \ $\ul{\beta}_i$ \
for \ $i>n$. \ Let \ $n\colon\Sigma_0\maps{1-1}\omega$ \ be such that for any \
$\psi\in\Sigma_0$, \  $\psi$ \ has support \ $n(\psi)$ \ and has no constants \ $\ux_i$ \ for \ $i\ge n(\psi)$. \ 
A variation of requirement R-5 (in the proof of Lemma~\ref{closedk>1}) is described below and allows  player
$\mathbf{I}$ to ``postpone'' his function moves. Finally, player $\mathbf{I}$ must prove
that he is constructing a model of  the theory \ $\T^\overn$ \ and that he is interpreting his constants correctly; this
involves commitments \ $m_j$ \ made at move \ $j$ \ as in the proof of Lemma~\ref{closedk>1}. 

The payoff of the game \ $G_x$ \ is essentially the same as the one described in the proof of Lemma
\ref{closedk>1} (see requirements R-1 to R-9), except there is no analogue for R-6 and requirement R-5 becomes:
\be
\item[R-5:] If \ $\sigma$ \ and \ $\tau$ \ are \ $\Sigma_0$ \ formulae of the form \ $(\iota v\phi(v)\in\F_n)$ \ and
\ $(\iota v\psi(v)\in\F_m)$, \ respectively, and if \ $\sigma,\tau$ \ are in \ $T^*$, \ then
\be
\item $m_{n(\sigma)}\ge n(\sigma)$ \ and \ $m_{n(\tau)}\ge n(\tau)$
\item $\left(\iota v\phi(v)\,\Eq^{s,t}\,\iota v\psi(v)\right)$  is in \ $T^*$ \ iff \ $\overcore\models
\f_{m_{n(\sigma)}}\,\Eq^{s,t}\,\f_{m_{n(\tau)}}$, \ 
for all \ $s\in(n+m)^n\!\uparrow$ \ and all \ $t\in(n+m)^m\!\uparrow$.
\ee
\ee
As before, requirement R-5 forces player $\mathbf{I}$ to prove that his model is mouse iterable. 
At move \ $n(\sigma)$, \ player $\mathbf{I}$ will play an integer \ $m_{n(\sigma)}$ \ predicting the move at which he will
present a function \ $\f_{m_{n(\sigma)}}$ \ fulfilling R-5. 

The definition of $x$--honesty is, in essence, as defined in the proof of Lemma~\ref{closedk>1} (see H-1 to H-8).
For instance, in the analogue for condition H-5, since player $\mathbf{I}$ can postpone his function moves (see the above
R-5), one must define \ $\tau=\{(\h_i,\f_{m_{n(\sigma_i)}}): i\le m\}$.
\setcounter{clm}{0}
\begin{clm}\label{claim1.2} The set \ $\{(x,u) : u \text{ is an $x$--honest position of length  $\n$}\}$ \ is \
$\boldface{\Sigma}{\omega}(\H^{\max\{i_0,\n\}})$.
\end{clm}
\begin{proof} The argument that all of the conditions of $x$--honesty are \
$\boldface{\Sigma}{\omega}(\H^{\max\{i_0,\n\}})$ \ proceeds as in the proof of Claim~\ref{claim1} of Lemma~\ref{closedk>1}. 
Although for  H-5,  Theorem~{\defontautwo} of \cite{Part1} is used to show that this condition is \
$\boldface{\Sigma}{\omega}(\H^{\max\{i_0,\n\}})$.
\end{proof}

\begin{clm}\label{claim2.2} For all \ $x\in\R$ \ and all \ $u$, \ the following are equivalent:
\be
\item $u$ \ is $x$--honest
\item $u$ \ is a winning position for player $\mathbf{I}$ in $G_x$.
\ee
\end{clm}
\begin{proof} The proof that (1) and (2) are equivalent is very similar to the proof of Claim~\ref{claim2} in
Lemma~\ref{closedk>1}. We note that in the proof of $(1)\Rightarrow(2)$, one takes advantage of player $\mathbf{I}$'s option to
postpone his function moves. Specifically, this option allows player $\mathbf{I}$ to choose \ $\f\in\F^{\overcore}$ \ so as to
ensure continued satisfaction of condition H-5 and, at the same time, obey the rule requiring his function moves \
$\f_j$ \ to be in \ $\F^{\H^j}$.
\end{proof}
Thus, the proof of Lemma~\ref{closedk=1limit} is finished.
\end{proof} 
\begin{lemma} If \ $\overcore\models\AD$ \ and \ $P\subseteq\R$ \ is \ $\boldface{\Sigma}{1}(\H^{i})$ \ for some \
$i\in\omega$, \ then \ $P$ \ has a scale which is \ $\boldface{\Sigma}{1}(\overcore)$.
\end{lemma}
\begin{proof} Suppose that \ $P\subseteq\R$ \ is \ $\boldface{\Sigma}{1}(\H^{i})$ \ for some  \ $i\in\omega$. \ It follows 
from the definition of \ $\H^{i+1}$ \ that for some  real \ $y$ \ there is a  \ $\Sigma_1$ \ formula in the parameters \
$y, \G, \rho^{}_{\overcore},p^{}_{\overcore},\beta_0, \beta_1, \dots, \beta_i,\beta_{i+1}$ \ that defines \ $P$ \ in \
$\H^{i+1}$. \ Lemma~\ref{closedk=1limit} now implies there is a scale on \ $P$ \ which is \
$\boldface{\Sigma}{1}(\overcore)$.
\end{proof}

This completes Case 2 in our proof of Theorem~\ref{goodcovering}. Our final case now follows.

{\bf Case 3: $\mathbf{k=1}$ and $\mathbf{\widehat{\OR}^{\overcore}}$ is a successor ordinal.} \ Since \
$\widehat{\OR}^{\overcore}$ \ is a successor ordinal, let \ $\nu$ \ be the ordinal such that  \
$\widehat{\OR}^{\overcore}=\nu+1$. \ Hence, \
${\OR}^{\overcore}=\omega\nu+\omega$. \ It follows that \ $n(\core)=0$.
\ To see this, suppose for a contradiction that \ $\overn=n(\core)\ge 1$. \ Thus, \
${\OR}^{\overcore}=\omega\rho^{}_{\core^{\overn-1}}$ \ and  Lemma~{\sigmacard} of \cite{Part1} 
asserts that \ ${\OR}^{\overcore}$ \ is a \ $\boldface{\Sigma}{1}(\core^{n-1})$--cardinal. 
However, \ ${\OR}^{\overcore}=\omega\nu+\omega$ \ is clearly not a \
$\boldface{\Sigma}{1}(\core^{\overn-1})$--cardinal. This contradiction implies that \ $n(\core)=0$.
\ Hence, \ $\overcore=\core$ \ and \ $\rho^{}_\core=1$. \ Therefore, \ $\core$ \ is a 1--mouse (see
Definition~{\onemouse} of \cite{Part1}). So, \ $\core$ \ is an iterable premouse of the form
\ $\core=(J^\mu_{\nu+1}(\R),\R,\kappa^\core,\mu)$.\footnote{Theorem 4.4 of \cite{Crcm} implies that $\Sigma_1(\core)$ has
the scale property but, because \ $\core$ \ is the core of a weak real mouse, we will be able to prove that \
$\boldface{\Sigma}{1}(\core)$ \ also has the scale property (see the proof of Theorem~\ref{newthm} in the next section).}
\ Consequently, when \ ${k=1}$ \ and \ ${\widehat{\OR}^{\overcore}}$ \ is a successor ordinal, the proof of Theorem
\ref{goodcovering}  is similar to the proof of Theorem 4.4 (for the successor case) in \cite{Crcm}. We shall now show how to
extend the proof of Theorem 4.4 in \cite[pp. 268-278]{Crcm} to handle this case.

Recalling Definition 1.4 and Lemma 1.6 of \cite{Crcm},  we see that \ $J^\mu_{\nu+1}(\R)=\bigcup\limits_{n \in \omega}S_{\nu+n}^{\core}({\R})$. \ Let \
$S_n =S_{\nu+n}^{\core}({\R})$ \ and let \ $\mathcal{S}_n=(S_n,\R,\kappa^\core,S_n\cap\mu)$. \ Recall that the domain of the structure \ $\core^\nu$ \ is \ $J^\mu_{\nu}(\R)=S_0$. \ Define functions \ $g_k$, \ inductively on \ $k$, \ by
\begin{align*}
g_0 &= f^{\core^{\,\nu}} \\
g_{k+1} &= G(g_k,S_k),
\end{align*}
where \ $f^{\core^\nu}$ \ is as in Corollary 1.8 \cite{Crcm} and \ $G$ \ is the rudimentary function given by Lemma 1.7 of
\cite{Crcm}. It follows that \ $g_k\colon [\omega\nu + k]^{<\omega}\times\R \maps{onto}S_k$ \ and since \ $G$ \ is rudimentary,
we can fix a recursive function \ $d\colon\omega \rightarrow \omega$ \ such that \ $g_k \in S_{d(k)}$ \ for each \
$k\in\omega$.

We define a canonical sequence \ $\seq{\H^i : i\in\omega}$ \ which will be a subsequence of \ $\seq{\mathcal{S}_n :
n\in\omega}$. \ At the same time we shall also define a sequence \ $\seq{\theta_i : i\in\omega}$ \ of \ $\Sigma_1$ \
formulae in \ $\Gamma$, \ a sequence \ $\seq{n_i : i\in\omega}$ \ of integers and a sequence \ $\seq{Q_{i} :
i\in\omega}$ \ of finite subsets of \ $\OR^{\core}$. 

First, define \ $Q_0=p^{}_{\core}$ \ and define \ $n_0$ \ to be the least integer such that \ $p^{}_{\core}\in S_{n_0}$. \ Let
\ $\H^0=\mathcal{S}_{n_0}$. \ We define \ $\theta_i$, \ $n_{i+1}$,  \ $Q_{i+1}$ \ and $\H^{i+1}$, \ by induction on \
$i$, \ as follows:

\smallskip
\noindent{\sl Definition of \ $\theta_i$\/}: \ Since \ $\H^i$ \ does not realize \ $\Gamma$, \ define \
$\theta_i$ \ to be the \ $<_\Gamma$\,--least \ $\Sigma_1$ \ formula \ $\theta(v)\in\Gamma$ \ such that 
\ $\H^i\not\models\theta(\seq{\G,w_0 })$.

\smallskip
\noindent{\sl Definition of \ $n_{i+1}$\/}: \ Because \ $\theta_i$ \ is \ $\Sigma_1$, \
let \ $\theta_i(v)=\exists u\psi_i(u,v)$ \ where \ $\psi_i$ \ is \ $\Sigma_0$. \
Since \ $\core\models\theta_i(\seq{\G,w_0})$, \ there is an \ $n\in\omega$ \
such that
\begin{equation}
\mathcal{S}_{n}\models\psi_i(g_k(Q,x),\seq{\G,w_0})\text{\ for some  
 and $x\in\R$}
\label{makeit2}\end{equation} 
for some \ $k\in\omega$ \ where \ $n>d(k)$ \ and some \ $Q\in[\omega\nu + k]^{<\omega}$.
\ Let \ $n_{i+1}$ \ be the least such \ $n>n_i$ \ satisfying (\ref{makeit2}).

\smallskip
\noindent{\sl Definition of \ $Q_{i+1}$\/}: \ Define \ $Q_{i+1}$ \ be the $\le_{\BK}$--least
\ $Q$ \ satisfying (\ref{makeit2}) with \ $n=n_{i+1}$.

\smallskip
\noindent{\sl Definition of \ $\H^{i+1}$\/}: \ Define \ $\H^{i+1}=\mathcal{S}_{n_{i+1}}$.

\smallskip
\begin{lemma}\label{lemma_k=1successor} Let \ $\seq{\H^i : i\in\omega}$, \ $\seq{\theta_i : i\in\omega}$ \ and \ $\seq{Q_i :
i\in\omega}$ \ be as above. Then \ $\bigcup\limits_{i\in\omega}\H^i=\core$ \ and for all \ $i\in\omega$
\be
\item[(a)] $\H^i\prec_{0}\H^{i+1}\prec_{0}\core$
\item[(b)] $\mathcal{H}^i$ \ is in \ $\boldface{\Sigma}{\omega}(\mathcal{H}^{i+1})$ 
\item[(c)] $\pow(\R)\cap\boldface{\Sigma}{\omega}(\H^i)\subseteq{C}$
\item[(d)] $\H^i\not\models\theta_i(\seq{\G,w_0})$ \ and \ $\H^{i+1}\models\theta_i(\seq{\G,w_0})$.
\ee
\end{lemma}
\begin{proof} For each \ $i\in\omega$, \ items (a)--(d) are clear.  The fact that \
$\H=\core$ \ is also clear.
\end{proof}
The sequence \ $\langle\H^i : i\in\omega\rangle$ \ is a suitable covering of \ $\core$. \ Let \ $\Lng_0 =
\{\,\en,\underline{\R},\underline{\kappa},\mu\,\}$.

\begin{lemma}\label{onto3} There is a partial \ $\boldface{\Sigma}{1}(\core)$ \ map of \ $\R$ \ onto \ $C$, \ the
domain of \ $\core$.
\end{lemma}
\begin{proof} Arguing as in the proof of Lemma~\ref{onto2}, one can show that the \ $\Sigma_1$ \ Skolem function \
$h_{\seq{\G,w_0}}$ \ (in the parameter \ $\seq{\G,w_0}$) \ maps \ $\R\times\R$ \ onto \ $C$.
\end{proof}

\begin{lemma}\label{closedk=1successor} Consider the parameters \ $\y, \G, Q_0,
Q_1, \dots, Q_\ell$ \ (as defined above) where \ $\y$ \ is a fixed real.  Let \ $\theta$ \ be a \ $\Sigma_1$ \ formula
in the language \ $\Lng_0$ \ augmented by these parameters. Fix \ $i_0\ge\ell$ \ and let \ $P\subseteq\R$ \ be defined by
\ $P(x)$ \ if and only if \ $\H^{i_0}\models \theta(x)$. \ Then \ $P$ \ has a closed game representation and if \ $\overcore\models\AD$, \ then there
is a \ $\boldface{\Sigma}{1}(\overcore)$ \ scale on \ $P$.
\end{lemma}
\begin{proof}[Sketch of Proof] The main ideas of the proof are very similar to those used in the proof
of the above Lemma~\ref{closedk>1} and in the proof of Theorem 4.4 in \cite{Crcm}. For this reason, we will just give an
outline of the argument. We want to construct a closed game representation of \ $P$ \ simple enough to yield a
\ $\boldface{\Sigma}{1}(\core)$ \ scale on \ $P$. \ Our closed game representation \ $x\mapsto G_x$ \
of \ $P$ \ will have the following property: If \ $P_n$ \ is defined by
\[P_n(x,u)\iff \textup{$u$ is a winning position for player $\mathbf{I}$ in $G_x$ of length $n$},\]
then \ $P_n$ \ is \ $\boldface{\Sigma}{\omega}(\H^{z(n)})$, \ for some fixed \ $z\colon\omega\to\omega$. \ Assuming \
$\overcore\models\AD$, \ Lemma~\ref{poffinmouse} and Lemma~\ref{lemma_k=1limit}(c) will then imply that the Moschovakis scale
\ $\seq{\varphi_i : i \in \omega}$ \ on \ $P$ \ exists. Furthermore, the prewellordering \ $\le_i$ \ of \ $\R$ \ induced by \
$\varphi_i$ \ is also \ $\boldface{\Sigma}{\omega}(\H^{z(i)})$. \ Lemma~\ref{lemma_k=1successor}(c) implies that \
$\le_i\,\in C$, the domain of \ $\core$, \ for all \ $i\in\omega$. \ Lemma~\ref{onto3} implies that any countable subset of \
$C$ \ is \ $\boldface{\Sigma}{1}(\core)$. \ Therefore, \ $\seq{\varphi_i : i \in\omega}$ \ is a \
$\boldface{\Sigma}{1}(\core)$ \ scale on \ $P$. 

So, to prove the lemma, it suffices to construct a closed game representation \ $x\mapsto G_x$ \ so that, for each \
$n\in\omega$, \ the predicate \ $P_n$ \ is first order definable over \ $\H^{k}$ \ for some \ $k\ge n$. \ Let \
$F^j=F^{\core}\cap S_j$ \ for each \ $j\in\omega$ \ (for a definition of $F^{\core}$, see section~\ref{sectwo} of this
paper).  In our construction of \ $G_x$ \ we want to force player $\mathbf{I}$ to describe the truth in \ $\core$ \ about \
$\G$  \ and the \ $Q_i$'s. \ To ensure that each \ $P_n$ \ is \ $\boldface{\Sigma}{\omega}(\H^{n})$, \ for some \ $k\ge n$, \
we impose some restrictions on player $\mathbf{I}$'s moves. First of all, player $\mathbf{I}$ must describe his model in the
language 
\[\Lng=\Lng_0\cup\{\uG\}\cup\{\ux_i,\ul{Q}_i,\ul{S}_i, \ul{g}_i: i\in\omega\}\]
where \ $\overn=n(\core)$. \ Let \ $\lh$ \ be a ``natural length'' function defined on all \ $\Sigma_{0}$ \ sentences \ in
the language \ $\Lng$ \ (for example, let \ $\lh(\psi)=\text{number of symbols in }\psi$). \ Also, let \
$\oz\colon\omega\to\omega$ \ be an increasing function whose properties will be described shortly. Let \ $T$ \ be the
theory consisting of the axioms listed on pages 270-271 in \cite{Crcm} (these axioms were used in the proof of Theorem 4.4
of \cite{Crcm})  minus axioms (7) and (13), together with the statement \ $(\uG\in[\OR]^{<\omega})$.  Player $\mathbf{I}$
must play a consistent and complete set of \ $\Sigma_{0}$ \  sentences in the language \ $\Lng$ \ extending the theory \
$T$.  \ At move
\ $i_0+1$ \ player $\mathbf{I}$ must play the \ $\Sigma_0$ \ sentence \ $\left((\ul{S}_{n_{i_0}},\mu)\models
\theta(\ux_0)\right)$. \  To prove that his model is an iterable premouse, player $\mathbf{I}$ must play functions from \
$F^j$ \ at each move \ $j$. \ A typical run of the game has the form
\begin{equation}
\begin{aligned}[c]
{}&{\mathbf{I}}\phantom{{\mathbf{I}}} \qquad T_0,s_0,f_0,m_0 \qquad \phantom{y_0} \qquad T_1,s_1,f_1,m_1 \qquad
\phantom{y_1}\quad
\\ {}&{\mathbf{I}}\mathbf{I} \qquad \phantom{T_0,s_0,f_0,m_0}\qquad y_0 \qquad \phantom{T_1,s_1,f_1,m_1} \qquad y_1 \quad 
\end{aligned}
\begin{gathered}[c]
{\cdots}
\end{gathered}
\tag*{$G_{x}$}\end{equation}
where 
\bi
\item $T_j$ \ is a finite set of $\Sigma_0$ sentences each of which has support \ $j$ \ and length \ $<\oz(j)$
\item $s_j\in\R^{<\omega}$, \ $f_j\in F^j$, \ $d(f_j)<j$, \ $j<m_j\in\omega$, \ and \ $y_j\in\R$. 
\ei
A \ $\Sigma_0$ \ formula \ $\varphi$ \ has {\it support\/} $n$ \ if it contains no constants \ $\ul{Q}_i,\ul{S}_i, \ul{g}_i$
\ for \ $i>n$. \ The function \ $d(f)$ \ is defined in section~\ref{sectwo}.  Let \ $n\colon\Sigma_0\maps{1-1}\omega$ \ be
such that for any \ $\psi\in\Sigma_0$, \ $\psi$ \ has support \ $n(\psi)$ \ and has no constants \ $\ux_i$ \ for \ $i\ge
n(\psi)$. \  An analogue of requirement R-5 (in the proof of Lemma~\ref{closedk>1}) is defined below and forces player
$\mathbf{I}$ to prove that his model is an iterable premouse. Finally, player $\mathbf{I}$ must interpret his constants
correctly; this involves commitments \ $m_j$ \ made at move \ $j$ \ as in the proof of Lemma~\ref{closedk>1}. 

The payoff of the game \ $G_x$ \ is essentially the same as the one described in the proof of Lemma
\ref{closedk>1} (see requirements R-1 to R-9), except there is no analogue for R-6 and requirement R-5 becomes the following
variation of requirement (2) in Definition 4.5 of \cite{Crcm}:
\be
\item[R-5:] If \ $\sigma$ \ and \ $\tau$ \ are \ $\Sigma_0$ \ formulae of the form \ $(\iota v\phi(v)\in F_n)^{\ul{S}_{i}}$ \
and
\ $(\iota v\psi(v)\in F_m)^{\ul{S}_{i'}}$, \ respectively, and if \ $\sigma,\tau$ \ are in \ $T^*$, \ then 
\[\text{$\left((\iota v\phi(v))^{\ul{S}_{i}}\,\le^{s,t}\,(\iota
v\psi(v))^{\ul{S}_{i'}}\right)^{\ul{\mathcal{S}}_{N}}$  is in
\ $T^*$ \ iff \
$\core\models
f_{n(\sigma)}\,\le^{s,t}\,f_{n(\tau)}$}\]
for all \ $s\in(n+m)^n\!\uparrow$ \ and all \ $t\in(n+m)^m\!\uparrow$, \ where \ $N=N(n, m, i, i')$.
\ee
The above R-5 forces player $\mathbf{I}$ to verify that his model is an iterable premouse (see \cite[Theorem 2.31]{Crcm}). The
relation \ $\le^{s,t}$ \ is define in \cite[p. 237]{Crcm}. \ Also, we let \
$\ul{\mathcal{S}}_{N}$ \ denote the term \ $(\ul{S}_N,\ul{S}_N\cap\mu)$. \ Finally, the recursive
function \ $N$ \ is described in \cite[remark on p. 270]{Crcm}. If \ $N = N(n, m, i, i')$, \ then \ $N$ \ is sufficiently large
so that the true structure \ ${\mathcal{S}}_N$ \ could legitimately decide the truth or falsehood of \ 
``$(\iota v\phi(v))^{{S}_{i}}\,\le^{s,t}\,(\iota v\psi(v))^{{S}_{i'}}$.''

The definition of $x$--honesty is obtained by combining the definition given in Lemma~\ref{closedk>1} (see H-1 to H-8)
together with the definition of honesty in \cite[Definiton 4.6]{Crcm}. For example, the analogue to condition H-5 becomes
\be
\item[H-5:] Let \ $\sigma_0,\dots, \sigma_i, \dots, \sigma_m$ \ enumerate those \ $\Sigma_0$ \ sentences that satisfy
the following three properties:  
\be
\item[(1)] $n(\sigma_i)<\n$
\item[(2)] $\sigma_i$ \ has the form \ $(\iota v\phi_i(v)\in F_{n_i})^{\ul{S}_{k_i}}$ \ where  \ $n_i,k_i\in\omega$
\item[(3)] $(\core, I_u)\models\sigma_i$. 
\ee 
For each \ $i\le m$ \ let \ $h_i\in F^{\core}_{n_i}$ \ be such that \ $(\core,
I_u)\models\phi_i(\h_i)^{\ul{S}_{k_i}}$. \ Let \ $\tau=\{(\h_i,f_{n(\sigma_i)}): i\le m\}$. \  Then there is \ $\Phi\colon
F^{\core}  \leadsto F^{\core}$ \ such that \ $\tau\subseteq\Phi$.
\ee
\begin{rmk} For a definition of \ $\Phi\colon F^{\core}  \leadsto F^{\core}$ \ in the above H-5,  \ 
see Definition 3.3 of \cite{Crcm}.
\end{rmk}
For another example, the analogue to condition H-7 becomes
\be
\item[H-7:] If \ $j<\n$ \ and the sentence \ $(\iota v\sigma(v)<_{\BK}\uG)^{\ul{\mathcal{S}}_{k}}$ \ is in \ $T_j$ \ for some \
$\Lng$--formula \ $\sigma$, \ then either 
\be
\item[(i)] $(\H^{m_j},I_u)\models\lnot\varphi(\seq{(\iota v\sigma(v))^{\ul{\mathcal{S}}_{k}},\ux_2})$ \ for some \ $\Pi_1$ \
formula \
$\varphi$ \ in the type \ $\Gamma$, \ or 
\item[(ii)] $(\H^{m_j},I_u)\models\lnot\varphi(\seq{(\iota v\sigma(v))^{\ul{\mathcal{S}}_{k}},\ux_2})$ \ for some \ $\Sigma_1$
\ formula
\ $\varphi$ \ which is one of the first \ $m_j$ \ elements in the type \ $\Gamma$.
\ee
\ee

\setcounter{clm}{0}
\begin{clm}\label{claim1.3} The set \ $\{(x,u) : u \text{ is an $x$--honest position of length  $\fn$}\}$ \ is \
$\boldface{\Sigma}{\omega}(\H^{k})$ \ for some \ $k\ge \fn$.
\end{clm}

\begin{proof}[Sketch of Proof] We shall discuss just the above conditions H-5 and H-7. The proof that H-5 is first order
over some \ $\H^{k}$ \ follows as in the proof of Theorem 4.4 in \cite{Crcm} for the successor case  (see \cite[pp.
277-8]{Crcm}). Now, because of our restrictions on the sentences in \ $T_j$ \ for \ $j<\fn$, \ we shall show that condition H-7 is first
order over some \ $\H^{k}$. \ Recall the sequences \ $\seq{\theta_i : i\in\omega}$ \ of \ $\Sigma_1$ \ formula in
\ $\Gamma$, \ $\seq{n_i : i\in\omega}$ \ defined in the construction of \ $\seq{\H^i : i\in\omega}$.
The real \ $\oz$ \ must be chosen so that for each \ $i\in\omega$
\[\lh(\psi_i(\ul{g}_a(\uQ_b,\ux_c),\seq{\uG,\ux_d}))<\oz(n_i)\]
for some (hence all) \ $a, b, c, d\in\omega$. \ Recall that \ $\H^{\fn}=\mathcal{S}_{n'}$ \ for some \ $n'$. \ Let \
$\Delta_{\fn}$ \ be the set of \ $\Sigma_0$ \ formulae, with one free variable, in the language \ $\Lng$ \ defined by
\[\Delta_{\fn}=\{ \sigma\in\Sigma_0 : \text{ $\sigma$ has support $<n'$ and $\lh(\sigma)<\fn$}.\}\]
A \ $\Sigma_0$ \ formula \ $\varphi$ \ with support support $<n'$ contains no constants \ $\ul{Q}_i,\ul{S}_i, \ul{g}_i$ \
for \ $i\ge n'$ \ and has no restrictions on the occurrence of any constants of the form \ $\ux_i$. \ We observe that for any
initial position \ $u$ \ of length \ $\fn$, \ all of the sentences in \ $\bigcup\limits_{j<\fn} T_j$ \ will be in \
$\Delta_{\fn}$. \ We also note that, modulo the constants \ $\ux_i$ \ and the variables of \ $\Lng$, \ the set \
$\Delta_{\fn}$ \ is finite. Given a real
\
$w$, \  define the interpretation \ $I_w$ \ in \ $\H^{\fn}$ \ as follows:
\[\begin{array}{l}
I_w(\ux_i)=w_i \text{ \ for all $i\in\omega$},\\
I_w(\uG)=\G,\\
I_w(\uQ_i)=Q_i \text{ \ for all $i<\fn$},\\
I_w(\ul{S}_i)=S_i \text{ \ for all $i<n'$},\\
I_w(\ul{g}_i)=g_i \text{ \ for all $i<n'$}\\
\end{array}\]
where \ $w_i$ \ denotes the $i$-th real encoded by \ $w$.\footnote{For example,
define
$w_i(j)=w(\seq{i,j})$ where  $\seq{i,j}$ is an integer recursively encoding the pair $(i,j)$.} 
Let \ $B_{\fn}$ \ be defined by
\[B_{\fn}=\{a\in H^{\fn} : (\H^{\fn},I_w)\models (a=\iota v \sigma(v)) \text{ for some } w\in\R \text{ and some } \sigma\in
\Delta_{\fn}\}.\] 
Because of the restrictions on the formulae in \ $\Delta_{\fn}$, \  there is a function \ $g$
\ in \ $\core$ \ such that \ $g\colon\R\maps{onto} B_{\fn}$. \ Consider the set \ $A$ \ defined by
\[A=\{\seq{K, \varphi, m} : K\in B_{\fn} \land \H^{\fn}\models(K<_{\BK}\G) \land \varphi\in\Gamma \land \H^m\models
\varphi(\seq{K,w_0})\}.\]
We shall show that \ $A$ \ is \ $\boldface{\Sigma}{\omega}(\H^k)$ \ for some \ $k$. \ Then one can
easily verify that condition H-7 is first order over \ $\H^k$. \ Using the above function \ $g$, \ it follows that there is a
function \ $f$ \ in \ $\core$ \ such that 
\[f\colon\R\maps{onto}\{K\in B_{\fn} : K<_{\BK}\G\}\times\Gamma.\]
Recall that \ $<_\Gamma$ \ is a fixed order of \ $\Gamma$ \ in order type \ $\omega$. \ For \ $\alpha,\beta\in\R$ \ define
\[\alpha\le^*\beta \iff f(\alpha)_0<_{\BK}f(\beta)_0 \lor (f(\alpha)_0=f(\beta)_0)\land f(\alpha)_1\le_{\Gamma}f(\beta)_1.\]
Since \ $f$ \ is in \ $\core$, \ we have that \ $\le^*$ \ is a prewellordering of \ $\R$ \ in \ $\core$. \ For an arbitrary 
\ $m\in\omega$, \ let \
\[A_m=\{\seq{K, \varphi} : \seq{K, \varphi, m}\in A\}.\]
Since \ $B_{\fn}$ \ is an element of \ $\core$, \ it follows (from the definition of \ $A$) \ that \ $A_m$ \ is in \ $\core$. \
Therefore, \ $f^{-1}(A_m)$ \ is in \ $\core$. \ The Coding Lemma (see \cite[see 7D.6]{Mosch}) implies that \ $f^{-1}(A_m)$
\ is \ $\boldfaceone{\Sigma}{1}(\le^*)$ \ for each \ $m\in\omega$. \ We remark that since \ $\le^*$ \ and \ $f^{-1}(A_m)$ \
are in \ $\core$, \ our assumption \ $\core\models\AD$ \ is sufficient in this case for applying the Coding Lemma. Because \
$\boldfaceone{\Sigma}{1}(\le^*)$ \ is closed under countable unions, we have that  \ $D\in\boldfaceone{\Sigma}{1}(\le^*)$ \
where
\[D=\bigcup\limits_{m\in\omega} f^{-1}(A_m)\times\{m\}.\]
Since \ $\core = \bigcup\limits_{n\in\omega}\H^n$ \ and
\[\seq{K, \varphi, m}\in A\iff (\exists y\in\R)(\seq{y,m}\in D \land f(y)_0=K\land f(y)=\varphi),\]
we conclude that \ $A$ is \ $\boldface{\Sigma}{\omega}(\H^{k})$ \ for some \ $k\ge \fn$.
\end{proof}

\begin{clm}\label{claim2.3} For all \ $x\in\R$ \ and all \ $u$, \ the following are equivalent:
\be
\item $u$ \ is $x$--honest
\item $u$ \ is a winning position for player $\mathbf{I}$ in $G_x$.
\ee
\end{clm}
\begin{proof} The proof that (1) and (2) are equivalent is very similar to the proof of Claim~\ref{claim2} in
Lemma~\ref{closedk>1}. The proof of the direction $(2)\Rightarrow(1)$ combines the proof of Theorem 4.7 of
\cite{Crcm} and the proof of this direction in  Lemma~\ref{closedk>1}.
\end{proof}
 
Claim~\ref{claim2.3} applied to the empty position implies that \ $x\mapsto G_x$ \ is a closed game representation of \ $P$. \
Claims~\ref{claim1.3} and \ref{claim2.3} \ imply, as stated at the beginning of the proof of Lemma~\ref{closedk=1successor},
that the resulting Moschovakis scale  on \ $P$ \ is \ $\boldface{\Sigma}{1}(\core)$. \ This
completes the proof of Lemma~\ref{closedk=1successor}.
\end{proof} 
\begin{lemma}\label{finally} If \ $\core\models\AD$ \ and \ $P\subseteq\R$ \ is \ $\boldface{\Sigma}{1}(\H^{i})$ \ for some \ $i\in\omega$,
\ then \ $P$ \ has a scale which is \ $\boldface{\Sigma}{1}(\core)$.
\end{lemma}
\begin{proof} Suppose that \ $P\subseteq\R$ \ is \ $\boldface{\Sigma}{1}(\H^{i})$ \ for some  \ $i\in\omega$. \ It follows 
that \ $P$ \ is in \ $C$, \ the domain of \ $\core$. \ Since the \ $\Sigma_1$ \ Skolem function \ $h_{\seq{\G,w_0}}$ \
(in the parameter \ $\seq{\G,w_0}$) \ maps \ $\R\times\R$ \ onto \ $C$ \ (see the proof of Lemma~\ref{onto3}), \ let \
$y',y''\in\R$ \ be so that \ $h_{\seq{\G,w_0}}(y',y'')=P$. \ Let \ $k\in\omega$ \ be sufficiently large so that this fact
holds in \ $\H^k$. \ Thus, \ $P$  \ is \ ${\Sigma_1}(\H^{k})$ \ in the parameters \ $y$ \ and \ $\G$ \ where \ $y=\seq{y',
y'', w_0}$ \ is a real which effectively encodes the triple \ $(y', y'',w_0)$. \ Lemma~\ref{closedk=1limit} 
implies there is a scale on \ $P$ \ which is \ $\boldface{\Sigma}{1}(\core)$.
\end{proof}

This completes the proof of Theorem~\ref{goodcovering}.
\end{proof}

\section{Weak real mice and scales}\label{finalsection}
We can now present a positive result on the existence of scales definable over a weak real premouse. 

\begin{theorem}[$\ZF+\DC$]\label{newthm}  Suppose that  \ $\mouse$ \  is a
weak real mouse satisfying \ $\AD$. \ Then \
$\boldface{\Sigma}{m}(\mouse)$ \ has the scale property, where \ $m=m(\mouse)$. 
\end{theorem}

\begin{proof} Let \ $\mouse$ \  be a
weak real mouse satisfying \ $\AD$. \ Let \ $\core=\core(\mouse)$, \ $\overn=n(\mouse)$ \ and \ $m=m(\mouse)$. \ 
Recall that \ $\overmouse=\mouse^\overn$ \ and \ $\overcore=\core^\overn$. \ Let \ $k\in\omega$ \ be such that \
$m=\overn+k$ \ and note that \ $k\ge 1$.  \ Because \ $\mouse$ \ is a mouse iterate of \ $\core$, \ Lemma 2.19 of \cite{Cfsrm}
implies that \ $\mouse$ \ and \ $\overcore$ \ have the same sets of reals. Thus, \ $\overcore\models\AD$. \ Suppose that \ $P$
\ is a \ $\boldface{\Sigma}{m}(\mouse)$ \ set of reals. Again, since  \ $\mouse$ \ is a mouse iterate of \ $\core$,  \ Lemma
2.19 of \cite{Cfsrm} asserts that \ $P$ \ is \ $\boldface{\Sigma}{k}(\overcore)$. \ Let  \
$\langle\mathcal{H}^i : i\in\omega\rangle$ \ be the suitable covering of \ $\overcore$ \ as defined in the proof of Theorem
\ref{goodcovering}.  Because \ $P$ \ is \ $\boldface{\Sigma}{k}(\overcore)$, \ let \ $a\in\overC$ \ (the domain of \
$\overcore$) and let \ $\varphi(u,v)$ \ be a \ $\Sigma_k$ \ formula in the language \ $\Lng_{\overn}$ \  such that
\[P(x) \iff \overcore\models\varphi(x,a)\]
for all \ $x\in\R$. \ Since \ $\overcore=\bigcup\limits_{i\in\omega}\H^i$, \ we shall assume without loss of generality that
\ $a\in H^0$ \ where \ $H^0$ \ is the domain of the structure \ $\H^0$. \ For each \ $i\in\omega$ \ define \ $P^i$ \ by
\[P^i(x) \iff \H^i\models\varphi(x,a)\]
for all \ $x\in\R$. \ It follows from the construction of the suitable covering that \ $P=\bigcup\limits_{i\in\omega}P^i$ \
(see the proof of Theorem~\ref{goodcovering}). For each \ $i$, \ let \ $\seq{\le^i_j : j\in\omega}$ \ be the \
$\boldface{\Sigma}{k}(\overcore)$ \  scale on \ $P^i$ \ constructed in the proof of Theorem~\ref{goodcovering}. We note  
that each \ $P^i$ \ and  \ $\le_j^i$  \ are elements in \ $\overcore$. \ Since Theorem~\ref{goodcovering} asserts that there
is a partial \ $\boldface{\Sigma}{k}(\overcore)$ \ map of \ $\R$ \ onto \ $\overC$, \ the map \ $(i,j)\mapsto(P^i,\le^i_j)$ \
is \ $\boldface{\Sigma}{k}(\overcore)$. \ Therefore, the union scale (see Definition~\ref{defunionscale}) is \
$\boldface{\Sigma}{k}(\overcore)$ \ on \ $P$ \ and Lemma 2.19 of \cite{Cfsrm} implies that this scale is also \
$\boldface{\Sigma}{m}(\mouse)$.
\end{proof}
\begin{corollary}[$\ZF+\DC$] Suppose that  \ $\mouse$ \  is a
weak real mouse satisfying \ $\AD$. Then every set of reals in \ $\mouse$ \  admits a \ $\boldface{\Sigma}{m}(\mouse)$ \
scale \ $\seq{ \le_i\, : i \in
\omega}$ \ each of whose norms \ $\le_i$ \ is in \ $\mouse$, \ where \ $m=m(\mouse)$.
\end{corollary}
\begin{proof} We shall assume the notation in the above proof of Theorem~\ref{newthm}. \ Let \ $P$ \ be a
set of reals in \ $\mouse$. \ Thus, \ $P$ \ is in \ $\overcore$. \ So,
\ $P$ \ is \ $\boldface{\Sigma}{k}(\H^i)$ \ for some \ $i\in\omega$. \ The proof of Theorem~\ref{goodcovering} implies that \
$P$ \ has a closed game representation and that the resulting Moschovakis scale \ $\seq{ \le_i\, : i \in
\omega}$ \ is as required.
\end{proof}

Recall that a weak real mouse is a proper initial segment of an iterable real premouse. Suppose that \ $\mouse$ \ is a real
mouse satisfying \ $\AD$, \ $m=m(\mouse)$ \ is defined, and \ $\mouse$ \ is not a proper initial segment
of an iterable real premouse. It would appear that one cannot apply Theorem~\ref{newthm}
to see if \ $\boldface{\Sigma}{m}(\mouse)$ \ has the scale property. The next corollary
asserts that, even in this case, one can apply Theorem~\ref{newthm}.

\begin{corollary}[$\ZF+\DC$] Let \ $\mouse$ \  be a
real mouse. Suppose that \ $m=m(\mouse)$ \ is defined and that \ $\mouse$ \ satisfies \ $\AD$. \ Let \ $\theta\in\OR$ \
be a multiple of \ $\omega^\omega$. \ If the mouse iterate \ $\mouse_\theta$ \ realizes a \ $\Sigma_m$ \ type not realized in
any proper initial segment of \ $\mouse_\theta$, \ then \ $\boldface{\Sigma}{m}(\mouse)$ \ has the scale property.
\end{corollary}
\begin{proof} Suppose that the mouse iterate \ $\mouse_\theta$ \ realizes a \ $\Sigma_m$ \ type not realized in any proper
initial segment of \ $\mouse_\theta$. \ By the proof of Theorem 2.49 of \cite{Cfsrm}, \ $\mouse_\theta$ \ is a proper
initial segment of an iterable real premouse. Therefore, \ $\mouse_\theta$ \ is weak. Theorem~\ref{newthm} asserts that \
$\boldface{\Sigma}{m}(\mouse_\theta)$ \ has the scale property. Lemma 2.19 of \cite{Cfsrm}
implies that \ $\boldface{\Sigma}{m}(\mouse)=\boldface{\Sigma}{m}(\mouse_\theta)$, \ as
pointclasses. We conclude that \ $\boldface{\Sigma}{m}(\mouse)$ \ has the scale property.
\end{proof}

When \ $\mouse$ \ is a weak real mouse and \ $m=m(\mouse)$, \ one can now make the observation that any set of reals
in \ $\boldface{\Sigma}{m}(\mouse)$ \ is the countable union of sets of reals in \ $\mouse$. \ The next two theorems follow
from the proof of Theorem~\ref{newthm}. 

\begin{theorem}\label{newcortwo} Suppose that  \ $\mouse$ \  is a
weak real mouse and let \ $m=m(\mouse)$. \ For any \ $\boldface{\Sigma}{m}(\mouse)$ \ set of reals \ 
$P$, \ there exists a total \ $\boldface{\Sigma}{m}(\mouse)$ \ 
map \ $k\colon\omega\to M$ \ such that \ 
$P=\bigcup\limits_{i\in\omega}k(i)$. 
\end{theorem}
\begin{theorem}\label{newcor1} Suppose that  \ $\mouse$ \  is a
weak real mouse and let \ $m=m(\mouse)$. \ Then for any set \ $P$ \ of reals, \
$P$ \ is \ in \ $\boldface{\Sigma}{m}(\mouse)$ \ if and only if \ $P=\bigcup\limits_{i\in\omega}k(i)$ \ for some map \
$k\colon\omega\to M$. 
\end{theorem}
\providecommand{\bysame}{\leavevmode\hbox to3em{\hrulefill}\thinspace}
\providecommand{\href}[2]{#2}

\end{document}